\documentclass[a4paper]{elsarticle}
\usepackage{lipsum}
\makeatletter
\def\ps@pprintTitle{%
 \let\@oddhead\@empty
 \let\@evenhead\@empty
 \def\@oddfoot{}%
 \let\@evenfoot\@oddfoot}
\makeatother
\usepackage{amssymb}
\usepackage{natbib}
\usepackage{amsmath}
\usepackage{cancel}
\usepackage{graphicx}
\usepackage[ansinew]{inputenc}
\usepackage{amsthm}
\usepackage[arrow, matrix, curve]{xy}
\usepackage{hyperref}

\newtheorem{thm}{Theorem}
\numberwithin{thm}{section}
\newtheorem{satz}[thm]{Theorem}
\newtheorem{Proposition}[thm]{Proposition}
\newdefinition{bemerkung}{Remark}
\numberwithin{bemerkung}{section}
\newtheorem{Lemma}[thm]{Lemma}

\newproof{pf}{Proof}
\newcommand{\ph}{\varphi}
\newcommand{\R}{\mathbb{R}}

\newcommand{\C}{\mathbb{C}}

\newcommand{\K}{\mathbb{K}}
\newcommand{\De}{\Delta_{p,q}}

\newcommand{\Pe}{\mathcal{P}}

\newcommand{\Se}{\mathcal{S}}
\newcommand{\g}{\mathfrak{g}}

\newcommand{\ort}{\mathfrak{o}}
\usepackage{amssymb}

\begin{document}

\begin{frontmatter}




\title{Towards a Classification of pseudo-Riemannian Geometries Admitting Twistor Spinors}


\author{Andree Lischewski}

\address{Department of Mathematics, Humboldt University, Rudower Chausse 25, 12489 Berlin, Germany}

\begin{abstract}
We show that for pseudo-Riemannian conformal structures a totally lightlike subspace fixed by the conformal holonomy representation is locally equivalent to having a Ricci-isotropic pseudo-Walker metric in the conformal class. This generalizes results obtained for lightlike lines and planes and naturally applies to parallel spin tractors resp. twistor spinors on conformal spin manifolds. In fact, it clarifies which twistor spinors are locally equivalent to parallel spinors. Moreover, we study the zero set of a twistor spinor using the curved orbit decomposition for parabolic geometries. Generalizing results from the Lorentzian case we can completely describe the local geometric structure of the zero set, construct a natural projective structure on it, and show that locally every twistor spinor with zero is equivalent to a parallel spinor off the zero set. An application of these results in low-dimensional split-signatures leads to a complete geometric description of local geometries admitting non-generic twistor spinors in signatures $(3,2)$ and $(3,3)$ which complements the well-known description of the generic case. In contrast to the generic case where generic geometric distributions play an important role, the underlying geometries in the non-generic case without zeroes turn out to admit integrable distributions.
\end{abstract}

\ead{lischews@math.hu-berlin.de}
\end{frontmatter}


\section{Introduction}
In this article we consider a space- and time-oriented, connected pseudo-Riemannian spin manifold of signature $(p,q)$. One can canonically associate to this setting  the real resp. complex spinor bundle $S^g$ with its Clifford multiplication, denoted by $\mu : TM \times S^g \rightarrow S^g$, and the Levi-Civita connection lifts to a covariant derivative $\nabla^{S^g}$ on this bundle. Besides the Dirac operator $D^g$, there is another conformally covariant differential operator acting on spinor fields, obtained by performing the spinor covariant derivative $\nabla^{S^g}$ followed by orthogonal projection onto the kernel of Clifford multiplication,
\[ P^g : \Gamma(S^g)  \stackrel{\nabla^{S^g}}{\rightarrow} \Gamma(T^*M \otimes S^g ) \stackrel{g}{\cong} \Gamma(TM \otimes S^g)  \stackrel{\text{proj}_{\text{ker}\mu}}{\rightarrow} \Gamma(\text{ker }\mu), \]
called the twistor operator. Elements of its kernel are called twistor spinors and they are equivalently characterized as solutions of the conformally covariant twistor equation
\[\nabla^{S^g}_X \ph + \frac{1}{n} X \cdot D^g \ph = 0 \text{    for all } X \in \mathfrak{X}(M). \]

In physics, twistor spinors appeared in the context of general relativity and were first introduced by R. Penrose in \cite{pen}. They became of interest in differential geometry   as T. Friedrich observed that special solutions of the twistor equation, the so called Killing spinors, are related to the first eigenvalue of the Dirac operator on a compact Riemannian spin manifold, see \cite{fr1}. Since then the twistor equation on Riemannian manifolds has been widely studied, e.g. in \cite{bfkg}. In particular, it is well-known that a Riemannian spin manifold admitting a twistor spinor without zeroes is conformally equivalent to an Einstein manifold which admits a parallel or a Killing spinor. The zero set in the Riemannian case has been widely studied (cf. \cite{ha, kr}). It consists of isolated points and if a zero exists, the spinor is conformally equivalent to a parallel spinor off the zero set. In contrast to the Riemannian and Lorentzian case (cf. \cite{lei, bl}), the investigation of the twistor equation in other signatures is widely open. The following general questions are of interest:
\begin{enumerate}
\item Which pseudo-Riemannian geometries admit nontrivial solutions of the twistor equation ?
\item How are further properties of twistor spinors related to the underlying geometries ? In particular, what are the possible shapes of the zero set $Z_{\ph} \subset M$ ?
\item How can one construct examples of manifolds admitting twistor spinors ?
\end{enumerate}
Recently, a $Spin^c$-version of the twistor equation became of interest in the context of the AdS/CFT correspondence in physics, see \cite{cks1, cks2}, and in this context the above questions help to distinguish pseudo-Riemannain manifolds on which supersymmetric conformal field theories can be placed. Obviously, the simplest subcase of twistor spinors are parallel spinors. \cite{kath} gives a complete classification of all non-locally symmetric, irreducible pseudo-Riemannian holonomy groups admitting parallel spinors. The other extremal case to irreducible acting holonomy is the case of a maximal holonomy invariant totally lightlike subspace. This leads to parallel pure spinors on pseudo-Riemannian manifolds which are studied in \cite{kath}. In split signatures, an explicit normal form of the metric is known. Furthermore, there are many examples and classification results for geometries admitting Killing spinors (cf. \cite{boh}, \cite{kath}). Another well-understood case are twistor spinors on Einstein spaces. \cite{bfkg} shows that in case of nonzero scalar curvature the spinor decoposes into a sum of two Killing spinors whereas in case of a Ricci-flat metric the spinor $D^g \ph$ is parallel. Also much is known about twistor spinors on Lorentzian manifolds. The most general result was obtained by F. Leitner in \cite{leihabil}. One can give up to conformal equivalence a complete list of local geometries admitting a twistor spinor off a certain singular set: Depending on the causal type of the associated conformal vector field $V_{\ph}$ one has a parallel spinor on a Brinkmann space, a local splitting into a Riemannian and Lorentzian factor, a Lorentzian-Einstein Sasaki structure or a Fefferman space. One can use this to deduce that the zero set of a twistor spinor with zero on a Lorentzian manifold consists either of isolated images of null-geodesics and off the zero set one has a parallel spinor on a Brinkmann space, or the zero set consists of isolated points and off the zero set one has a local splitting $(\R,-dt^2) \times (N,h)$, where the last factor is Riemannian Ricci-flat Kaehler, in the conformal class.\\
\cite{hs} indicates that in signatures higher than Lorentzian there are new interesting relations between twistor spinors and constructions of conformal structures out of projective structres. However, as there is no complete classification of manifolds admitting twistor spinors, one often restricts oneself to small dimensions in order to find out which geometries play a role there. \cite{br} classifies metrics admitting parallel spinor fields in small dimensions. It is moreover known that a Riemannian 3-manifold admitting a twistor spinor is conformally flat, and a Riemannian 4-manifold with twistor spinor is selfdual (\cite{bfkg}). In Lorentzian geometry, there is a classification of all local geometries admitting twistor spinors without zeroes and constant causal type of the associated conformal vector field $V_{\ph}$ for dimensions $n \leq 7$, which can be found in \cite{lei} or \cite{bl}. In signature $(2,2)$, anti-self-dual four manifolds with paralel real spinor have been studied in \cite{dus}. Furthermore, \cite{hs} presents a Fefferman construction which starts with a 2-dimensional projective structure and produces geometries carrying two pure spin tractors with nontrivial pairing which leads to $Hol(M,c) \subset SL(3,\R) \subset SO^+(3,3)$. 
\cite{hs1} investigates (real) \textit{generic} twistor spinors in signature $(3,2)$ and $(3,3)$, being twistor spinors satisfying additionally that the constant (!) $\langle \ph , D \ph \rangle \neq 0$ (signature $(3,3)$ is also discussed in \cite{br2}). They are shown to be in tight relationship to so called generic 2-distributions on 5-manifolds resp. generic 3-distributions on 6-manifolds, that means every generic twistor spinor gives rise to a generic distribution, and conversely, given a manifold with generic distribution, one can canonically construct a conformal structure admitting a twistor spinor, and these two constructions are inverse to each other.\\
\newline
Twistor spinors are objects of conformal geometry and (except the Riemannian case) all mentioned results in the pseudo-Riemannian context are established by making use of conformal tractor calculus and by equivalently describing twistor spinors as parallel sections in the spin tractor bundle associated to a conformal spin manifold as presented in \cite{baju} or \cite{lei}. In this setting, geometries admitting twistor spinors are equivalently characterized as those conformal spaces $(M,c)$ where the lift of the conformal holonomy group $Hol(M,c) \subset SO(p+1,q+1)$ to $Spin(p+1,q+1)$ stabilizes a nontrivial spinor. A problem closely related to the twistor equation is therefore the classification of pseudo-Riemannian conformal holonomy groups which is completely solved only in the Riemannian case (cf. \cite{baju}). In arbitrary signatures, one knows a conformal analogue of the local de-Rham/Wu-splitting theorem (cf. \cite{leihabil}) and all holonomy groups acting transitive and irreducible on the Moebius sphere were classified in \cite{alt}. The most involved case is the situation when the holonomy representation fixes a totally lightlike subspace $H \subset \R^{p+1,q+1}$. The associated local geometries are only known in cases dim $H \leq 2$ (\cite{baju, ln}).

In this article we study precisely this classification problem for conformal holonomy groups, i.e. the case that a totally lightlike subspace of dimension $\geq 2$ is fixed by the holonomy representation and show in Proposition \ref{ct} that if on a conformal manifold $(M,c)$ there exists a totally lightlike, $k$-dimensional parallel distribution in the standard tractor bundle, then every point of some open and dense subset admits a neighborhood $U$, a metric $g \in c_U$ and a $k-1$-dimensional totally lightlike distribution $L \subset TU$ such that
\begin{align}
Ric^g(TU) \subset L,\label{33} \\ 
L \text{ is parallel wrt. }\nabla^g. \label{44}
\end{align}
Conversely, if $U \subset M$ is an open set equipped with a metric $g \in c_U$ and a $k-1$-dimensional totally lightlike distribution $L \subset TU$ such that (\ref{33}) and (\ref{44}) hold, then $L$ gives rise to a $k-$dimensional totally lightlike, parallel distribution in the standard tractor bundle over $U$.\\
\newline
In the rest of the article we apply this result to the classification problem for twistor spinors, and we show that a large class of twistor spinors is locally equivalent to parallel spinors off a certain singular set. In fact, if we have a parallel spin tractor, we can associate via the holonomy principle a holonomy invariant spinor $v \in \Delta_{p+1,q+1}$ (up to conjugation). As an application of Proposition 3.2 we show that if the kernel $H_v$ of this spinor under Clifford multiplication with vectors from $\R^{p+1,q+1}$ is nontrivial, the associated twistor spinor is locally conformally equivalent to a parallel spinor off a singular set (Proposition \ref{tms}). Also the converse is true. In the remainder of the article we then present two main applications of these results: First, we are able to clarify the local structure of the zero set of a twistor spinor in arbitrary signatures: In Theorem \ref{zss} we show that for $\ph \in \Gamma(S^g)$a twistor spinor with zero, the zero set $Z_{\ph}$ is an embedded, totally lightlike submanifold. Moreover, for every $x \in Z_{\ph}$ there are open neighborhoods $U$ of $x$ in $M$ and $V$ of $0$ in $T_xM$ such that
\begin{align*}
Z_{\ph} \cap U = \text{exp}_x \left(\text{ker }D^g \ph(x) \cap V \right). 
\end{align*}
Besides, we show that the conformal structure canonically induces a torsion-free projective structure on the zero set of a twistor spinor (Proposition \ref{lps}). In this regad we mention \cite{der} where a similar statement is proved for zero set components of certain conformal vector fields.  As a second application we study (real) twistor spinors in small dimensions. We are able to classify geometries admitting non-generic twistor spinors in signature $(3,2)$ and $(3,3)$ which complements the analysis of the generic case from \cite{hs}. We prove in Proposition \ref{propos} that real twistor half-spinors in signature $(2,2)$ without zeroes and real twistor (half-)spinors without zeroes in signatures $(3,2)$ and $(3,3)$ satisfying that $\langle \ph , D^g \ph \rangle \equiv 0$ are locally conformally equivalent to parallel spinors (off a singular set). Their associated distributions ker $\ph \subset TM$ are integrable (off a singular set). Finally, we can also obtain some results in the less studied signatures $(4,2)$ and $(4,3)$.\\
\newline
This article is organized as follows:
In the first section we  provide the elaboration of the fundamental principles and methods to work on the classification problem for twistor spinors in pseudo-Riemannian geometry, including algebraic preliminaries and basic facts about conformal tractor calculus and conformal holonomy. In section 3 we prove the classification result for conformal structures admitting totally lightlike, holonomy-invariant subspaces. Applications to twistor spinors with zeroes are then presented in section 4 where we make use of the so called curved orbit decomposition for arbitrary parabolic geometries from \cite{cgh}. Finally, in section 5 we discuss twistor spinors in low dimensions. Algebraic observations regarding the orbit structure of the spinor module $\Delta_{p,q}$ in low dimensions as known from \cite{br} directly relate the general previous results to concrete statementes in low dimensions.

\section{Twistor spinors and associated objects}
\subsection{The real and complex spinor module}
We consider $\R^{p,q}$, that is, $\R^n$, where $n=p+q$, equipped with a scalar product $\langle \cdot, \cdot \rangle_{p,q}$ of index $p$, given by $\langle e_i, e_j \rangle_{p,q} = \epsilon_i \delta_{ij}$, where $(e_1,...,e_n)$ denotes the standard basis of $\R^n$ and $\epsilon_i \in \{\pm 1 \}$ are fixed. In general, we should think of $\langle \cdot, \cdot \rangle_{p,q}$ as being the pseudo-Euclidean standard scalar product of index $p$, i.e. $\epsilon_i=-1$ for $1 \leq i \leq p$ and $\epsilon_i=+1$ for $p+1 \leq i \leq n$. However, in order to simplify the following calculations, we shall work with this more general notion of $\R^{p,q}$. We denote by $Cl_{p,q}$ the Clifford algebra of $(\R^{n},- \langle \cdot, \cdot \rangle_{p,q})$ and by $Cl_{p,q}^{\C}$ its complexification. It is the associative real or complex algebra with unit multiplicatively generated by $(e_1,...,e_n)$ with the relations \[ e_ie_j+e_je_i=-2 \langle e_i,e_j \rangle_{p,q}. \] 
It is well-known (see \cite{lm, har}) that if $p-q \not \equiv 1 $mod $4$, there is (up to equivalence) exactly one irreducible real representation of $Cl_{p,q}$. If $p-q \equiv 1 $\text{mod} $4$, there are precisely two inequivalent real irreducible representations of $Cl_{p,q}$. Furthermore, ${Cl}^{\C}_{p,q}$ admits up to equivalence exactly one irreducible complex representation in case $n$ is even and two such representations if $n$ is odd. In case that there are two equivalence classes of irreducible real or complex representations, they can be distinguished by the unit volume element as presented in \cite{lm}: Let $\omega_{\R}:= e_1 \cdot....\cdot e_{n}  \in Cl_{p,q}$ and
$\omega_{\C}:= (-i)^{\left[\frac{n+1}{2}\right]-p} \omega_{\R} \in Cl^{\C}_{p,q}$. If $p-q \equiv 1$ mod $4$, each irreducible real representation of $Cl_{p,q}$ or $Cl_{p,q}^{\C}$ maps $\omega_{\R}$ to $Id$ or $-Id$. Both possibilities can occur and the resulting representations are inequivalent. The analogous statements are true in the complex case for $Cl_{p,q}^{\C}$ and $n$ odd (cf. \cite{ba81}).
This opens a way to distinguish a up to equivalence unique real resp. complex irreducible representation for all Clifford algebras $Cl_{p,q}$ and $Cl_{p,q}^{\C}$ by requiring that $\omega$ is mapped to $Id$ in case $n$ even ($\K=\C$) or $p-q \equiv 1$ mod $4$ ($\K=\R$).

\begin{bemerkung} \label{ds}
For concrete calculations we shall make use of the following irreducible, complex representation of $Cl^{\C}_{p,q}$:
Let $E,T,g_1$ and $g_2$ denote the $2 \times 2$ matrices
\begin{align*} 
E = \begin{pmatrix} 1 & 0 \\ 0 & 1 \end{pmatrix} \text{ , } T = \begin{pmatrix} -1 & 0 \\ 0 & 1 \end{pmatrix} \text{ , } g_1 = \begin{pmatrix} 0 & i \\ i & 0 \end{pmatrix} \text{ , } g_2 = \begin{pmatrix} 0 & -1 \\ 1 & 0 \end{pmatrix}.
\end{align*}
Furthermore, let \[ \tau_j =\begin{cases}  1 &  \epsilon_j = 1, \\ i & \epsilon_j = -1. \end{cases} \]
Let $n=2m$. In this case, ${Cl}^{\C}(p,q) \cong M_{2^m}(\C)$ as complex algebras, and an explicit realisation of this isomorphism is given by
\begin{align*}
\Phi_{p,q} (e_{2j-1})&= \tau_{2j-1} \cdot E \otimes...\otimes E \otimes g_1 \otimes \underbrace{T \otimes...\otimes T}_{(j-1) \times},\\
\Phi_{p,q} (e_{2j})  &= \tau_{2j} \cdot E \otimes...\otimes E \otimes g_2 \otimes \underbrace{T \otimes...\otimes T}_{(j-1) \times}.
\end{align*}
Let $n=2m+1$. In this case, there is an isomorphism $\widetilde{\Phi}_{p,q} : {Cl}^{\C}({p,q}) \rightarrow M_{2^m}(\C) \oplus M_{2^m}(\C)$, given by
\begin{align*}
\widetilde{\Phi}_{p,q} (e_j) &= (\Phi_{p,q-1}(e_j),\Phi_{p,q-1}(e_j)) \text{,  } j=1,...,2m ,\\
\widetilde{\Phi}_{p,q} (e_{2m+1}) & = \tau_{2m+1} (iT \otimes...\otimes T, -iT \otimes...\otimes T),
\end{align*}
and $\Phi_{p,q}:=pr_1 \circ \widetilde{\Phi}_{p,q}$ is an irreducible representation mapping $\omega_{\C}$ to $Id$.
\end{bemerkung}

Fixing an irreducible real or complex representation $\rho:Cl^{(\C)}_{p,q} \rightarrow End(\De)$ and restricting it to the spin group $Spin(p,q) \subset Cl_{p,q} \subset Cl_{p,q}^{\C}$ yields a representation of $Spin^{(+)}(p,q)$\footnote{By $Spin^{(+)}(p,q)$ we denote $Spin(p,q)$ or it's identity component $Spin^+(p,q)$.} on the space of real or complex spinors $\De \in \{\Delta_{p,q}^{\R}, \Delta_{p,q}^{\C} \}$, called the \textit{real or complex spinor representation}. One possible realisation in the complex case is $\De^{\C}=\C^{2^m}$, where $n=2m+1$ or $n=2m$ (cf. Remark \ref{ds}). In case $n$ even $(\K=\C)$ or $p-q \equiv 0$ mod $8$ ($\K=\R$), $\De$ splits into the sum of two inequivalent $Spin(p,q)$ representations $\De^{\pm}$ according to the $\pm 1$ eigenspaces of $\omega$ (cf. \cite{har, ba81}). In our realisation from Remark \ref{ds} one can find these half spinor modules as follows:
Let us denote by $u(1)$ the vector $\begin{pmatrix}1 \\ 0 \end{pmatrix} \in \C^2$, by $u(-1)$ the vector $\begin{pmatrix} 0 \\ 1 \end{pmatrix} \in \C^2$ and set $u(\epsilon_1,...,\epsilon_m):=u(\epsilon_m) \otimes...\otimes u(\epsilon_1)$ for $\epsilon_{\nu} = \pm1$. Then we have
\begin{align*}
\Delta^{\C,\pm}_{p,q} = \text{span} \{ u(\epsilon_1,...,\epsilon_m) \mid \prod_{\nu = 1}^m \epsilon_{\nu} = \pm 1 \}.
\end{align*}
Note further that $Cl_{p,q}^{(\C)}$ acts on $\De$ via the representation $\rho$, and as $\R^n \subset Cl_{p,q} \subset Cl^{\C}_{p,q}$, this defines the Clifford multiplication $ X \cdot \ph := \rho(X)(\ph)$ of a vector by a spinor which naturally extends to a multiplication by $k$-forms: Letting 
$\omega = \sum_{1 \leq i_1 <...< i_k \leq n} \omega_{i_1...i_k} e^{\flat}_{i_1} \wedge...\wedge e^{\flat}_{i_k} \in \Lambda^k_{p,q}:=\Lambda^k \left(\R^{p,q}\right)^*$ and $\ph \in \De$, we set
\begin{align*} \omega \cdot \ph := \sum_{1 \leq i_1 <...< i_k \leq n} \omega_{i_1...i_k} e_{i_1} \cdot...\cdot  e_{i_k} \cdot \ph \in \De.  \end{align*}
In the \textit{split signatures} $(m+1,m)$ and $(m,m)$ the spinor representations $\De^{\C}$ are real (cf. \cite{har}). More precisely, there exists a real structure $\alpha$ on $\De^{\C}$ commuting with Clifford multiplication such that the real spinor module can be realised to be $\De^{\R} = \{ \ph \in \De^{\C} \mid \alpha(\ph) = \ph \} \subset \De^{\C}$. Setting $\epsilon_j = (-1)^j$ in these cases, one sees that our realisation of complex Clifford multiplication from Remark \ref{ds} is then given by real matrices in the split signatures. Therefore, it restricts to a real action of $Cl_{p,q}$ on $\R^{2^m}$ and the real structure $\alpha$ is simply given by complex conjugation.\\
\newline
Next, we will define $Spin^+(p,q)$-invariant inner products on $\De^{\C}$, following \cite{ba81}. To this end, we choose an irreducible representation of $Cl^{\C}_{p,q}$ such that $\De^{\C}$ can be realised to be $\C^{2^{[n/2]}}$. Let $( \cdot, \cdot )_{\De}$ denote the standard scalar product on this space. Then the bilinear form $\langle \cdot, \cdot \rangle_{\De^{\C}}$\footnote{If $q>p$ we instead work with a scalar product which involves $e_{i_j}$ with $\epsilon_{i_j}=1$ in the definition. It has analogous properties, cf. \cite{har}.}, given by
\begin{align}
\langle u, v \rangle_{\De^{\C}} = d \cdot (e_{i_1} \cdot...\cdot e_{i_p} \cdot u, v )_{\De^{\C}}, \label{2}
\end{align}
where $d$ is some power of $i$ depending on $p,q$ and the concrete realisation of the representation only, $i_1<...<i_p$ and $\epsilon_{i_1}=...=\epsilon_{i_p}=-1$, is a Hermitian scalar product on $\De^{\C}$. If $p,q >0$, it has neutral signature and it holds that 
\begin{align} \langle X \cdot u, v \rangle_{\De^{\C}} + (-1)^p \langle u, X \cdot v \rangle_{\De^{\C}} = 0. \label{fg}\end{align} for all $u,v \in \De^{\C}$ and $X \in \R^n$.
In the real case, we can proceed analogous (cf. \cite{sp}) by choosing an irreducible real representation of $Cl_{p,q}$ such that as vector space $\De^{\R}$ can be realised to be $\R^N$ for some $N$ (cf. \cite{har}). We then let $( \cdot, \cdot )_{\De^{\R}}$ denote the standard scalar product on this space and define $\langle \cdot, \cdot \rangle_{\De^{\R}}$ as in (\ref{2}), where me may now set $d=1$. (\ref{fg}) still holds in the real case. Moreover, $\langle \cdot , \cdot \rangle_{\De^{\R}}$ is symmetric if $p=0,1$ mod $4$ with neutral signature ($p \neq0$ and $q \neq 0$) or it is definite ($p=0$ or $q=0$). In case $p=2,3$ mod $4$, the pair $(\De^{\R},\langle \cdot , \cdot \rangle_{\De^{\R}})$ is a symplectic vector space.\\
\newline
There is an important decomposition of $\Delta_{p+1,q+1}$ into $Spin(p,q)-$modules. Let $(e_0,...,e_{n+1})$ denote the standard basis of $\R^{p+1,q+1}$. We introduce lightlike directions $e_{\pm} := \frac{1}{\sqrt{2}}(e_{n+1} \pm e_0)$.  One then has a decomposition $\R^{p+1,q+1} = \R e_- \oplus \R^{p,q}\oplus \R e_+ $ of $\R^{p+1,q+1}$ into $O(p,q)-$modules. We define the annihilation spaces $Ann(e_{\pm}):=\{ v \in \Delta_{p+1,q+1} \mid e_{\pm}\cdot v = 0 \}$. It follows that for every $v \in \Delta_{p+1,q+1}$ there is a unique $w \in \Delta_{p+1,q+1}$ such that $v=e_- w + e_+ w$, leading to a decomposition
\begin{align}
\Delta_{p+1,q+1} = Ann(e_-) \oplus Ann(e_+). \label{fs}
\end{align}
As $x e_{\pm}=-e_{\pm}x$ for all $x \in \R^{p,q} \cong \text{span}(e_1,...,e_n) \subset \R^{p+1,q+1}$ we see that $\R^{p,q}$ and $Spin(p,q)$ act on Ann$(e_{\pm})$. We can thus realise $\Delta_{p,q}$ as being $Ann(e_{\pm})$. Now fix an isomorphism $\alpha: Ann(e_-) \rightarrow \Delta_{p,q}$ of $Spin(p,q)$-representations. Then there is an induced isomorphism $\beta: Ann(e_+) \rightarrow \Delta_{p,q}$, $v \mapsto \alpha(e_- v)$ and an isomorphism
\begin{align}
{\Delta_{p+1,q+1}}_{|Spin(p,q)} & \cong \Delta_{p,q} \oplus \Delta_{p,q}, \\
v=e_+w+e_-w & \mapsto (\alpha(e_-w),\alpha(e_-e_+w)) \label{deco}
\end{align}
of $Spin(p,q)$ modules. One calculates that wrt. this decomposition the scalar product $\langle \cdot , \cdot \rangle_{\Delta_{p+1,q+1}}$ is given by
\begin{align}
\langle \begin{pmatrix} v_1 \\ w_1 \end{pmatrix}, \begin{pmatrix} v_2 \\ w_2 \end{pmatrix} \rangle_{\Delta_{p+1,q+1}} = -\frac{\delta^p}{\sqrt{2}} \left( \langle v_1,w_2 \rangle_{\Delta_{p,q}} + (-1)^p \langle w_1,v_2 \rangle_{\Delta_{p,q}}\right) \label{zz}
\end{align}
where $v_j,w_j \in \Delta_{p,q}$ for $j=1,2$ and $\delta=i$ in case $\mathbb{K}=\mathbb{C}$ and $\delta=1$ in the real case.\\
\newline
In general, the orbit structure of $\De$ under the $Spin^+(p,q)$ action becomes very complicated as $n=p+q$ increases. However, in small dimensions the orbits are well understood (cf. \cite{br}), and there is one distinguished orbit which turns out to be of particular importance here, namely the so called \textit{pure spinors} (cf. \cite{br, kath, pure}). In order to define them, we follow \cite{kath}. Let $m=[n/2]$. First, we consider the complex case: The Clifford mutliplication can be extended to a complex bilinear map $\C^n \times \De^{\C} \rightarrow \De^{\C}$. To each spinor $v \in \De^{\C}$ we associate the subspaces 
\begin{align*} \text{ker}_{\C} v := \{ X \in \C^n \mid X \cdot v = 0 \} \text{ and } \text{ker } v := \{ X \in \R^n \mid X \cdot v = 0 \}. \end{align*}
One checks that ker$_{\C} v$ is isotropic with respect to the complex linear extension $\langle \cdot, \cdot \rangle_{p,q}^{\C}$ of $\langle \cdot, \cdot \rangle_{p,q}$, and in particular, ker $v$ is isotropic with respect to $\langle \cdot, \cdot \rangle_{p,q}$. Clearly, dim ker $v$ is an $Spin(p,q)$-orbit invariant and it holds that $\lambda(g)(\text{ker }v) \subset \text{ker }v$ for all $g \in Spin(p,q)$, where $\lambda: Spin(p,q) \rightarrow SO(p,q)$ denotes the double covering map.
A complex spinor $v \in \De^{\C}$ is said to be \textit{pure} if dim$_{\C} \text{ ker}_{\C}v = \left\lceil \frac{n}{2}\right\rceil$, i.e., if its kernel under (extended) Clifford multiplication is a maximally isotropic subspace. In this case, dim$_{\R} \text{ker } v$ is called the \textit{real index} of $v$.
Next, we consider the real case and pure spinors in $\De^{\R}$. The notion of a real pure spinor can be developed for all signatures $(p,q)$ as explained in \cite{br}, but we are only interested in pure spinors in the \textit{split} signatures $(m,m)$ and $(m+1,m)$, and in these cases pure spinors in $\De^{\R}$ can be defined using the complex definition (cf. \cite{kath}): Consider the split signatures $(m,m)$ and $(m+1,m)$ and the inclusion of the real spinor module $\De^{\R} \subset \De^{\C} = \De^{\R} \oplus i \De^{\R}$ (using a real structure as explained before). Then a spinor $v \in \De^{\R}$ is called \textit{(real) pure} if it is the real or imaginary part of a pure spinor in the complexified module $\De^{\C}$ which has real index $m$.
So in the following, when talking about pure spinors, we mean either the complex case or real pure spinors in split signature. If $n=2m$ the set of pure spinors in $\De$ forms precisely one orbit under the $Spin^+(p,q)$ action, whereas in case $n=2m$ pure spinors form one orbit in each half spinor module $\De^{\pm}$.\\
Given a real pure spinor $\chi \in \De^{\R}$ in split signature, \cite{kath} shows that up to conjugation in $Spin^+(p,q)$ the stabilizer of $\chi$ under the $Spin^+(p,q)$ action is given by
\begin{align} 
Stab_{\chi}Spin^+(p,q) = SL(m) \ltimes N,\label{stab}
\end{align}
where $N$ is a certain nilpotent group.

\subsection{Associated forms to a spinor}

In the Lorentzian case one can associate to every nonzero spinor a nonvanishing vector, the so called Dirac current. Generalizing this construction, we associate to every spinor $\chi \in \De$ a series of forms $\alpha_{\chi}^k \in \Lambda^k_{p,q}$, $k \in \mathbb{N}$, so called \textit{algebraic Dirac forms}, given by 
\begin{align}
\langle \alpha_{\chi}^k,\alpha \rangle_{p,q} := d_{k,p} \cdot \langle \alpha \cdot \chi, \chi \rangle_{\Delta_{p,q}} \textit{  } \forall \alpha \in \Lambda^k_{p,q}. \label{6}
\end{align}
$d_{k,p}$ is a nonzero constant depending on the chosen representation but not depending on $\chi$, ensuring that the so defined form is indeed a real form.
The following properties of these forms turn out to be important and are easily checked:
\begin{Proposition} \label{10}
Let $\chi \in \De$ and $k \in \mathbb{N}$.
\begin{enumerate}
\item $\alpha^p_{\chi}=0 \Leftrightarrow \chi = 0$
\item  $\alpha_{\chi}^k = d_{k,p}  \sum_{1 \leq i_1 < i_2 <...<i_k \leq n} \epsilon_{i_1}...\epsilon_{i_k} \langle e_{i_1}\cdot...e_{i_k}\cdot \chi, \chi  \rangle_{\Delta_{p,q}} e^{\flat}_{i_1} \wedge...\wedge e^{\flat}_{i_k}$ for some constant $d_{k,p}$.
\item Equivariance: $\alpha^k_{g \cdot \chi } = \lambda(g) (\alpha^k_{\chi})$ for all $k \in \mathbb{N}$, $g \in Spin^+(p,q)$ and $\chi \in \Delta_{p,q},$ where $\lambda:Spin(p,q) \rightarrow SO(p,q)$ denotes the double covering map.
\end{enumerate}
\end{Proposition}
There is an important relation between the structure of $\alpha_{\chi}^p$ and $\text{ker } \chi$:

\begin{Lemma}\label{0}
Let $\chi \in \Delta_{p,q}\backslash \{0 \}$ and let $k:= \text{dim }\text{ker }\chi( \leq p)$. Then $\alpha_{\chi}^p$ can be written as
\begin{align}
\alpha_{\chi}^p = l_1^{\flat} \wedge ... \wedge l_k^{\flat} \wedge \widetilde{\alpha}, \label{maxi}
\end{align}
where $l_j \in \R^{p,q}$ for $1 \leq j \leq k$ such that $\text{span }\{l_1,...,l_k\} = \text{ker }\chi$ (in particular, this implies that the $l_j$ are lightlike and mutually orthogonal), $\widetilde{\alpha} \in \Lambda^{p-k}\left(\left(\text{ker }\chi\right)^{\bot}\right)$ and (\ref{maxi}) is \textit{maximal} in the sense that there exists no lightlike vector $l_{k+1}$ being orthogonal to $l_i$ for $1 \leq i \leq k$ such that $\alpha_{\chi}^p = l_1^{\flat} \wedge ...\wedge l_k^{\flat} \wedge l_{k+1}^{\flat} \wedge \widetilde{\widetilde{\alpha}}$.
Moreover, whenever $\alpha_{\chi}^p$ can be written as in (\ref{maxi}) for mutually orthogonal lightlike vectors $l_1,...,l_k$, it follows that $l_1,...,l_k \in \text{ker }\chi$.
\end{Lemma}
The Proof of Lemma \ref{0} is a straightforward generalization of the proof for the case $p=1$ as presented in \cite{lei}. Lemma \ref{0} generalize well known facts about the associated Dirac current $V_{\chi}$ to a spinor $\chi \in \Delta_{1,n-1}$ in the Lorentzian case from \cite{lei}: It holds that $||V_{\chi}||^2=0$ implies that $V_{\chi} \cdot \chi = 0$ being is a special case of Lemma \ref{0}, and $V_{\chi}$ is always causal.

\begin{bemerkung} \label{rem}
All possible algebraic Dirac forms $\alpha_{\ph}^2$ for $0 \neq \ph \in \Delta^{\C}_{2,n-2}$ have been classified in \cite{leihabil}. Precisely one of the following cases occurs:
\begin{enumerate}
\item $\alpha_{\ph}^2 = l_1^{\flat} \wedge l_2^{\flat}$, where $l_1,l_2$ span a totally lightlike plane in $\R^{2,n-2}$.
\item $\alpha_{\ph}^2 = l^{\flat} \wedge t^{\flat}$ where $l$ is lightlike, $t$ is a orthogonal timelike vector.
\item $\alpha_{\ph}^2 =\omega_0$ (up to conjugation in $SO(2,n-2)$), where $\omega_0$ is the standard Kaehler form on $\R^{2,n-2}$. In this case $n=2m$ and $Stab_{\alpha_{\ph}^2} O(2,n-2) \subset U(1,m-1)$.
\item There is a nontrivial Euclidean subspace $E \subset \R^{2,n-2}$ such that ${\alpha_{\ph}^2}_{|E} = 0$ and $\alpha_{\ph}^2$ is the standard Kaehler form on the orthogonal complement $E^{\bot}$ of signature $(2,2m)$ (again, this is up to conjugation in $SO(2,n-2)$). In this case $Stab_{\alpha_{\ph}^2} O(2,n-2) \subset U(1,m) \times O(n-2(m+1))$.
\end{enumerate}
Lemma \ref{0} implies that the first case occurs iff ker $ \ph$ is maximal, i.e. 2-dimensional. The second case occurs iff this kernel is one-dimensional whereas the last two cases can only occur if the kernel under Clifford multiplication is trivial.
\end{bemerkung}

\subsection{The twistor equation on spinors} \label{tes}
We fix some notations about basic objects from spin geometry (following \cite{ba81}) and recall the definition and properties of twistor spinors as introduced in \cite{bfkg}. Let $(M,g)$ be a space-and time-oriented, connected pseudo-Riemannian spin manifold of  index $p$ and dimension $n=p+q \geq 3$. By $\Pe^g_+$ we denote the $SO^+(p,q)$-principal bundle of all space-and time-oriented pseudo-orthonormal frames\footnote{In the following, these frames are simply referred to as
 pseudo-orthonormal.} $s=(s_1,...,s_n)$. A spin structure  of $(M,g)$ is then given by a $\lambda-$reduction $(\mathcal{Q}^g_+,f^g)$ of $\Pe^g_+$ to  the group $Spin^+(p,q)$. In particular, $f^g:\mathcal{Q}^g_+ \rightarrow \Pe^g_+$ is a 2-fold covering. The associated bundle $S^g:=\mathcal{Q}^g_+ \times_{Spin^+(p,q)} \De$ is called the \textit{real or complex spinor bundle}. In case that $\De=\De^+ \oplus \De^-$, it holds that $S^g=S^{g,+} \oplus S^{g,-}$, and one then has the notion of half-spinor fields. The algebraic objects introduced in the last section define fibrewise Clifford multiplication $\mu : T^*M \otimes S^g \rightarrow S^g$ and an inner product $\langle \cdot , \cdot \rangle_{S^g}$. Clearly, the properties of $\langle \cdot , \cdot \rangle_{\De}$ translate into corresponding properties of $\langle \cdot , \cdot \rangle_{S^g}$. Finally, the Levi Civita connection $\nabla^g$ on $(M,g)$, considered as a bundle connection $\omega^g \in \Omega^1(\Pe_+^g,\ort(p,q))$, lifts to a connection $\widetilde{\omega}^g \in \Omega^1(\mathcal{Q}^g_+,\mathfrak{spin}(p,q))$ which in turn induces a covariant derivative $\nabla^{S^g}$ on $S^g$. Locally, $\nabla^{S^g}$ is given by the formula 
\begin{align*}
\nabla^{S^g}_X \ph = X(\ph) + \frac{1}{2} \sum_{1 \leq k < l \leq n} \epsilon_i \epsilon_j g(\nabla^g_X s_k,s_l) s_ks_l \cdot \ph, 
\end{align*}
for $\ph \in \Gamma(S^g)$ and $X \in \mathfrak{X}(M)$, where $s=(s_1,...,s_n)$ is any oriented local pseudo-orthonormal frame. The composition of $\nabla^{S^g}$ with Clifford multiplication defines the Dirac operator
\[ D^g : \Gamma(S^g) \stackrel{\nabla^{S^g}}{\rightarrow} \Gamma(T^*M \otimes S^g ) \stackrel{g}{\cong} \Gamma(TM \otimes S^g) \stackrel{\mu}{\rightarrow} \Gamma(S^g),\]
whereas performing the spinor covariant derivative $\nabla^{S^g}$ followed by orthogonal projection onto the kernel of Clifford multiplication gives rise to the \textit{twistor operator} $P^g$
\[ P^g : \Gamma(S^g)  \stackrel{\nabla^{S^g}}{\rightarrow} \Gamma(T^*M \otimes S^g ) \stackrel{g}{\cong} \Gamma(TM \otimes S^g)  \stackrel{\text{proj}_{\text{ker}\mu}}{\rightarrow} \Gamma(\text{ker} \mu). \]
Spinor fields $\ph \in \text{ker }P^g$ are called \textit{twistor spinors} and they are equivalently characterized as solutions of the \textit{twistor equation}
\[\nabla^{S^g}_X \ph + \frac{1}{n} X \cdot D^g \ph = 0 \text{    for all } X \in \mathfrak{X}(M). \]
The twistor operator is conformally covariant: Letting $\widetilde{g}=e^{2 \sigma}g$ be a conformal change of the metric, it holds that $P^{\widetilde{g}} \widetilde{\varphi} = e^{-\frac{\sigma}{2}} \left( P^g(e^{-\frac{\sigma}{2}} \varphi) \right)\widetilde{}$, where $\widetilde{}$ denotes the natural identification of $\mathcal{Q}^g$ with $\mathcal{Q}^{\widetilde{g}}$. In particular, $\varphi \in \Gamma(S^g)$ is a twistor spinor with respect to $g$ if and only if the rescaled spinor $e^\frac{\sigma}{2} \widetilde{\varphi} \in \Gamma(S^{\widetilde{g}})$ is a twistor spinor with respect to $\widetilde{g}$. Moreover, the dimension of the space of twistor spinors is conformally invariant and bounded by $2^{[n/2]+1}$.
If $\ph \in \Gamma(S^g)$ is a twistor spinor, it holds that 
\begin{align}
\nabla^{S^g}_X D^g \ph = \frac{n}{2}K^g(X)\cdot \ph. \label{shou}
\end{align}
where $K^g=\frac{1}{n-2} \left(\frac{scal^g}{2(n-1)}g-Ric^g \right)$ denotes the Shouten tensor.\\
\newline
The previous observations show that twistor spinors are in fact objects of conformal geometry. One is therefore intended to develop a concept describing twistor spinors if one has only given a conformal class $c=[g]$ instead of a single metric $g \in c$. One elegant approach to do this, is making use of the conformal tractor calculus as presented in \cite{baju} or \cite{feh}. To this end, let $(M,c)$ be a connected, space-and time-oriented conformal manifold of signature $(p,q)$ and dimension $n=p+q \geq 3$.
We call a frame $(s_1,...,s_n)$ over $x \in M$ a \textit{conformal frame} if there is $g \in c$ such that the vectors $s_1,...,s_n$ are pseudo-orthonormal with respect to this metric. Collecting all these frames, we obtain the \textit{conformal frame bundle}
$(\Pe_+^0 , \pi^0 , M ; CO^+(p,q) )$ with structure group the identity component of the conformal group $CO(p,q)\cong \R^+ \times O(p,q)$.
Using the general theory of parabolic geometries (cf. \cite{cs}), one shows that the oriented conformal structure $(M,c)$ is \textit{equivalently} encoded in a normal parabolic geometry\footnote{In fact, for our purposes one does not need to introduce the general concept of parabolic geometries for this equivalence. For an explicit construction via first prolongation we refer to \cite{baju} or \cite{feh}} $(\Pe_+^1,\pi, M,\omega^{nc})$ of type $(G,P)$ over $M$, where we have the following objects: $G=SO(p+1,q+1)$, and the parabolic subgroup $P \subset G$ is defined as follows: Let $(e_0,...,e_{n+1})$ denote the standard basis of $\R^{p+1,q+1}$, introduce two lightlike directions by setting $e_{\pm}:=\frac{1}{\sqrt{2}}(e_{n+1} \pm e_0)$ and let $P:=Stab_{\R^+e_-}G$, where $G$ acts on $\R^{p+1,q+1}$ via the standard matrix action. These algebraic objects describe the flat model $(G \rightarrow G/P, \omega^{MC})$ (with $\omega^{MC}$ being the Maurer Cartan form) for conformal structures, being a double cover $\widehat{Q}^{p,q}:=G/P$ of the (pseudo-)Moebius sphere $Q^{p,q}$ equipped with a flat conformal structure $\widehat{c}$. We set $\widehat{C}^{p,q}:=(\widehat{Q}^{p,q},\widehat{c})$. For a general conformal structure $(\Pe_+^1,\pi, M,\omega^{nc})$ the Cartan connection $\omega^{nc} \in \Omega^1(\Pe^1_+,\g)$ on $\Pe^1_+$ is uniquely determined by certain normalisation conditions (cf. \cite{baju}) and called the \textit{normal conformal Cartan connection}. It describes the deviation from the flat model. It´s extension to a principal bundle connection on the principal $G$-bundle $\overline{\Pe^1_+}:=\Pe_+^1 \times_P G$ induces a covariant derivative $\nabla^{nc}$ on the standard tractor bundle $\mathcal{T}(M):=\Pe^1_+ \times_P \R^{p+1,q+1}$. Furthermore, $\langle \cdot, \cdot \rangle_{p+1,q+1}$ induces a bundle metric on $\mathcal{T}(M)$, and $\nabla^{nc}$ turns out to be metric. Therefore, it seems natural to define the \textit{conformal holonomy} of the conformal structure to be the holonomy of this connection \footnote{For a slightly different, but in our case equivalent approach to define the holonomy of a Cartan connection we refer to \cite{baju}}:
\[ Hol_x(M,c):=Hol_x(\mathcal{T}(M),\nabla^{nc}) \subset SO^+(\mathcal{T}(M)_x,\langle \cdot, \cdot \rangle) \cong SO^+(p+1,q+1). \]
The null line $I = \mathbb{R} e_- \subset \R^{p+1,q+1}$ which defines the parabolic subgroup $P$ induces a filtration $I \subset I^{\bot} \subset \R^{p+1,q+1}$ which leads to a filtration $\mathcal{I} \subset \mathcal{I}^{\bot} \subset \mathcal{T}(M) $ of $\mathcal{T}(M)$.\\
Fixing a metric $g \in c$ leads to a natural reduction $\sigma^g: \Pe^g \rightarrow \Pe^1$ to the structure group $SO^+(p,q)$. The standard representation of $SO^+(p+1,q+1)$ on $\R^{p+1,q+1}$ splits into the direct sum of 3 $SO^+(p,q)$ representations:
\begin{align*}
\R^{p+1,q+1} \cong \R \oplus \R^{p,q} \oplus \R \text{ with }ae_- + Y + b e_+ \mapsto (a,Y,b).
\end{align*}
This in turn leads to an isomorphism
\begin{align}
\mathcal{T}(M) \stackrel{g}{\cong} \underline{M} \oplus TM \oplus \underline{M}. \label{ide}
\end{align}
With respect to this identification we have that $\mathcal{I} = \underline{M} \oplus 0 \oplus 0$ and  $\mathcal{I}^{\bot} = \underline{M} \oplus TM \oplus 0$. In particular, we can use $g$ to identify sections $s \in \Gamma(\mathcal{T}(M))$ with triples $(\alpha, Y, \beta)$, where $\alpha, \beta \in C^{\infty}(M)$ and $Y \in \mathfrak{X}(M)$. Under this identification, the bundle metric is given by
 \begin{align} \langle (\alpha_1, Y_1, \beta_1), (\alpha_2, Y_2, \beta_2) \rangle_{\mathcal{T}(M)} = \alpha_1 \beta_2 + \alpha_2 \beta_1 + g(Y_1,Y_2), \label{cc1} \end{align}
 and one has the following formulas for the metric description of the tractor connection $\nabla^{nc}$ and its curvature ${R}^{\nabla^{nc}}$:
 \begin{align} \nabla_X^{nc} \begin{pmatrix} \alpha \\ Y \\ \beta \end{pmatrix} = \begin{pmatrix} X(\alpha) + K^g(X,Y) \\ \nabla_X^g Y + \alpha X - \beta K^g(X)^{\sharp} \\ X(\beta) - g(X,Y) \end{pmatrix} \text{, }R^{\nabla^{nc}}_{X_1,X_2}\begin{pmatrix} \alpha \\ Y \\ \beta \end{pmatrix} = \begin{pmatrix} C^g(X_1,X_2)Y \\ W^g(X_1,X_2)Y - \beta C^g(X_1,X_2)^{\sharp} \\0 \end{pmatrix}, \label{cc2} \end{align}
where $C^g(X,Y):=\nabla^g_X(K^g)(Y)-\nabla^g_Y(K^g)(X)$ defines the Cotton-tensor and $W^g$ is the Weyl-tensor of the conformal structure.
Under a conformal change $\widetilde{g} = e^{2 \sigma}g$, the metric representation of a standard tractor transforms according to (cf. \cite{baju})
\begin{align}
\begin{pmatrix} \alpha \\ Y \\ \beta \end{pmatrix} \mapsto \begin{pmatrix} \widetilde{\alpha} \\ \widetilde{Y} \\ \widetilde{\beta} \end{pmatrix}= \begin{pmatrix} e^{- \sigma} (\alpha - Y(\sigma) - \frac{1}{2}\beta ||\text{grad}^g \sigma ||^2_g \\ e^{- \sigma} (Y + \beta \text{grad}^g \sigma) \\ e^{\sigma} \beta \end{pmatrix}. \label{tra}
\end{align}

An analogous first prolongation procedure can be carried out in the conformal \textit{spin} setting (cf. \cite{baju}). Let $CSpin^{(+)}(p,q)=\R^+ \times Spin^{(+)}(p,q)$ denote the (identity component of) the conformal spin group, coming together with a double covering $\lambda^0:CSpin^+(p,q) \rightarrow CO^+(p,q)$. 
$(\widetilde{G} = \lambda^{-1}(G), \widetilde{P}=\lambda^{-1}(P))$ is the pair on which conformal spin structures are modelled as parabolic geometries. The Cartan geometry $(\widetilde{G} \rightarrow \widetilde{G} / \widetilde{P} \cong \widehat{Q}^{p,q}, \omega^{MC})$ is the flat model and can be viewed as the space $\widehat{C}^{p,q}$ equipped with a conformal spin structure (cf. \cite{leihabil}). For concrete calculations we use a realisation of the flat model preseted in \cite{lei}:  $\widehat{Q}^{p,q}$ is isomorphic to the set of time-oriented null directions in $\R^{p+1,q+1}$. It is naturally embedded in $\R^{p+1,q+1}$ via
\begin{align}
i :\widehat{Q}^{p,q} \hookrightarrow \R^{p+1,q+1} \text{, where }\R_+ \cdot x \mapsto \sqrt{\frac{2}{\langle x,x \rangle_{n+2}}} \cdot x \label{dat}
\end{align} 
and where $\langle \cdot,\cdot \rangle_{n+2}$ denotes the standard Euclidean inner product. One checks that $i(\widehat{Q}^{p,q}) = S^p \times S^q \subset \R^{0,p+1} \times \R^{0,q+1}$. It holds that $\widehat{c}=[i^*\langle \cdot, \cdot \rangle_{p+1,q+1}]$, yielding the conformally flat conformal spin manifold $\widehat{C}^{p,q}=(\widehat{Q}^{p,q},\widehat{c})$,  which realises the flat model for conformal spin structures of index $p$. Suppose now that $(M,c)$ carries a conformal spin structure, being a $\lambda^0$-reduction $(\mathcal{Q}^0_+,f^0)$ of the bundle $\Pe^0_+$ to $CSpin^+(p,q)$. In analogy with the previous case, this geometric structure is via first prolongation equivalently encoded in a parabolic geometry $(\mathcal{Q}^1_+,\widetilde{\pi},M,\widetilde{\omega}^{nc})$ of type $(\widetilde{G},\widetilde{P})$ such that one has the following double coverings:
\begin{align*}
\begin{array}{llllll}
(\mathcal{Q}^g_+;Spin^+(p,q))& \hookrightarrow &(\mathcal{Q}^0_+;CSpin^+(p,q))& \leftrightarrow &(\mathcal{Q}^1_+,\widetilde{\omega^{nc}})&(\widetilde{G},\widetilde{P})\\
 f^g \downarrow & & f^0 \downarrow && f^1 \downarrow & \lambda \downarrow  \\
(\Pe^g_+;SO^+(p,q))& \hookrightarrow & (\Pe^0_+;CO^+(p,q))& \leftrightarrow &(\Pe^1_+,\omega^{nc})&(G,P)\\
\end{array}
\end{align*}
The normal conformal spin connection $\widetilde{\omega}^{nc} \in \Omega^1(\mathcal{Q}^1_+,\mathfrak{spin}(p+1,q+1))$ induces a covariant derivative - also denoted by $\nabla^{nc}$ - on the (real or complex conformal) \textit{spin tractor bundle} $\mathcal{S}_{\mathcal{T}}(M):= \mathcal{Q}^1_+ \times_{\widetilde{P}} \Delta_{p+1,q+1}$. Furthermore, one has in analogy to the metric setting an inner product $\langle \cdot, \cdot \rangle_{\mathcal{S}}$ on this bundle, and a pointwise Clifford multiplication $\mu(X,\psi):=X \cdot \psi$ of sections $X \in \Gamma(\mathcal{T}(M))$ and spinor fields $\psi \in \Gamma( \mathcal{S}_{\mathcal{T}}(M))$. \\
Fixing a metric $g \in c$ leads to a reduction $\sigma^g: \mathcal{Q}_+^g \rightarrow \mathcal{Q}^1_+$ of $(\mathcal{Q}^1_+,\widetilde{P})$ to $(\mathcal{Q}^g_+,Spin^+(p,q))$. We let $\overline{\mathcal{Q}^1_+}$ denote the enlarged $Spin^+(p+1,q+1)$-principal bundle, and as $\mathcal{S}_{\mathcal{T}}(M) \cong \overline{\mathcal{Q}^1_+} \times_{Spin^+(p+1,q+1)} \Delta_{p+1,q+1}$, we may use $g$ to identify
\begin{align*}
Q^g_+ \times_{\rho \circ i_{cs}} \Delta_{p+1,q+1} & \cong \mathcal{S}_{\mathcal{T}}(M), \\
[l,v] & \mapsto [[\sigma^g(l),e],v],
\end{align*}
where $i_{cs}:Spin(p,q) \hookrightarrow Spin(p+1,q+1)$ denotes the natural inclusion, and $e \in Spin(p+1,q+1)$ is the neutral element. The decomposition (\ref{deco}) of $\Delta_{p+1,q+1}$ induces projections $\text{proj}^g_{\pm}: \mathcal{S}_{\mathcal{T}}(M) \rightarrow \overline{\mathcal{Q}^1_+} \times_{Spin(p,q)} Ann(e_{\pm})$ and an vector bundle isomorphism
\begin{align}
\Phi^g: \mathcal{S}_{\mathcal{T}}(M) & \rightarrow S^g(M) \oplus S^g(M), \label{gdg}\\ 
[[\sigma^g(l),e],e_-w+e_+w] & \mapsto [l,\beta(e_+w)] + [l,\alpha(e_-w)]. 
\end{align}
One calculates that under this identification, $\nabla^{nc}$ is given by the expression (cf. \cite{baju})
\begin{align*}
\nabla^{nc}_X \begin{pmatrix} \ph  \\ \phi \end{pmatrix} = \begin{pmatrix} \nabla_X^{S^g} & -X \cdot \\ \frac{1}{2}K^g(X) \cdot & \nabla^{S^g}_X \end{pmatrix} \begin{pmatrix} \ph  \\ \phi \end{pmatrix}.
\end{align*}
Together with (\ref{shou}) this yields a reinterpretation of twistor spinors in terms of conformal Cartan geometry:
\begin{satz} \label{90}
Let $(M,c)$ be a connected, space- and time-oriented conformal spin manifold of dimension $n\geq 3$. 
For any metric $g \in c$, the vector spaces of twistor spinors in $\Gamma(S^g)$ and parallel sections in $\Gamma (\mathcal{S}_{\mathcal{T}}(M))$ are naturally isomorphic via
\begin{align*}
 \text{ker }P^g \rightarrow \Gamma(S^g(M) \oplus S^g(M)) \stackrel{\left(\Phi^g\right)^{-1}}{\cong } \Gamma(\mathcal{S}_{\mathcal{T}}(M))\text{,   } 
 \ph \mapsto \begin{pmatrix} \ph \\  -\frac{1}{n}D^g \ph \end{pmatrix} \stackrel{\left(\Phi^g\right)^{-1}}{\mapsto} \psi \in Par(\mathcal{S}_{\mathcal{T}}(M), \nabla^{nc}), 
 \end{align*}
 i.e. a spin tractor $\psi \in \Gamma(\mathcal{S}_{\mathcal{T}}(M))$ is parallel iff for one (and hence for all $g \in c$), it holds that $\ph:= \Phi^g \circ \text{proj}_+^g \psi \in \text{ker }P^g$ and $D^g \ph = -n \cdot \Phi^g \circ \text{proj}_-^g \psi$.
\end{satz}
 
\subsection{The twistor equation on forms}
In the Lorentzian case, it has paid off to associate to every spinor field $\ph \in \Gamma(S^g)$ a vector field $V_{\ph} \in \mathfrak{X}(M)$, the so called Dirac current, as done in \cite{bl,lei}. The zero sets of these objects coincide, i.e. $Z_{\ph}=Z_{V_{\ph}}$. If $\ph$ is a twistor spinor, then $V_{\ph}$ is a conformal vector field and Lorentzian geometries admitting twistor spinors can partially be classified by studying the behaviour of $V_{\ph}$, cf. \cite{bl} for details. This procedure can be generalized to arbitrary signatures by making use of the nc-Killing form theory as presented in \cite{nc} or \cite{leihabil}. We list some later needed facts:\\
A global version of (\ref{6}) associates to a spinor field $\ph \in \Gamma(S^g)$ a  $k-$form $\alpha^k_{\ph} \in \Omega^k(M)$ for each $k \in \mathbb{N}$ with $Z_{\ph}=Z_{\alpha^p_{\ph}}$. In the special case $k=1$, it holds that $\alpha_{\ph}^1=V_{\ph}^{\flat}$. If $\ph$ is a twistor spinor, the forms $\alpha^k_{\ph}$ turn out to be \textit{normal conformal (nc-)Killing k-forms}, meaning that they are conformal Killing forms,
\begin{align*}
 \nabla^{g}_X \alpha_{\ph}^k - \frac{1}{k+1} \iota_X d \alpha_{\ph}^k + \frac{1}{n-k+1} X^{\flat} \wedge d^* \alpha_{\ph}^k &= 0,\text{ for all }X \in \mathfrak{X}(M),
\end{align*}
which satisfy additional normalisation conditions as to be found in \cite{lei}. One checks that if $\alpha=\alpha^k_{\ph}$ is a nc-Killing $k$-form wrt. $g$, then the rescaled form \begin{align} \widetilde{\alpha}:=e^{(k+1)\sigma}\alpha^k_{\ph}=\alpha^k_{e^{\frac{\sigma}{2}}\widetilde{\ph}} \label{8} \end{align} is a nc-Killing $k$-form wrt. $\widetilde{g}=e^{2\sigma}g$.\\
On the other hand, if we view the twistor spinor as parallel spin tractor $\psi \in Par(\Se,\nabla^{nc})$, we can also associated to this object a series of forms. In order to define them, we introduce the \textit{tractor k-form bundle} $\Lambda^k_{\mathcal{T}} (M):=\mathcal{P}^1_+ \times_{P} \Lambda^k_{p+1,q+1}$. Sections, i.e. elements of $\Omega^k_{\mathcal{P}^1}(M) :=\Gamma(\Lambda^k_{\mathcal{T}} (M))$ are called  \textit{tractor $k$-forms on $M$}. Clearly, the standard scalar product on $\Lambda^k_{p+1,q+1}$ induces a bundle metric on this space and the normal conformal Cartan connection $\omega^{nc}$ leads to a covariant derivative $\nabla^{nc}$. Again, (\ref{6}) can be applied pointwise, and in this way, a series of tractor forms $\alpha^k_{\psi}$ is associated to every spin tractor $\psi \in \Gamma(\Se)$. In the special case of $\psi$ being parallel, $\alpha^k_{\psi}$ turns out to be parallel as well. Parallel tractor $k-$forms are called \textit{(normal) twistor k-forms}.\\
Fixing a metric in the conformal class leads to the following description of tractor $k$-forms: First, note that every form $\alpha \in \Lambda_{p+1,q+1}^{k+1}$ decomposes into $\alpha = e_+^{\flat} \wedge \alpha_- + \alpha_0 + e^{\flat}_- \wedge e_+^{\flat} \wedge \alpha_{\mp} + e_-^{\flat} \wedge \alpha_+ $
for uniquely determined forms $\alpha_-,\alpha_+ \in \Lambda^k_{p,q}, \alpha_0 \in \Lambda^{k+1}_{p,q}$ and $\alpha_{\mp } \in \Lambda^{k-1}_{p,q}$. 
Whence, we can reduce $\Pe_+^1$ to $\Pe^g_+$ with structure group $SO^+(p,q)$ and see that there is an isomorphism $\Lambda^{k+1}_{\mathcal{T}} (M) \stackrel{g}{\cong} \Lambda^{k}(M) \oplus \Lambda^{k+1}(M) \oplus \Lambda^{k-1}(M) \oplus \Lambda^{k}(M)$,
and consequently, each tractor $(k+1)$-form $\alpha \in \Omega^{k+1}_{\mathcal{P}^1}(M)$ uniquely corresponds via a fixed metric $g \in c$ to a set of differential forms $\alpha \stackrel{g}{\leftrightarrow} (\alpha_- , \alpha_0 , \alpha_{\mp}, \alpha_+)$, where $\alpha_-,\alpha_+ \in \Omega^k(M), \alpha_0 \in \Omega^{k+1}(M), \alpha_{\mp} \in \Omega^{k-1}(M)$. We can also write this as
\begin{align} 
\alpha \stackrel{g}{=}s_-^{\flat} \wedge \alpha_- + \alpha_0 + s_-^{\flat} \wedge s_+^{\flat} \wedge \alpha_{\mp}+ s_+^{\flat} \wedge \alpha_+, \label{ag}
\end{align}
i.e. the $s_{\pm}$ are global lightlike sections in the line bundles in $\mathcal{T}(M) \stackrel{g}{\cong}\underline{M} \oplus TM \oplus \underline{M}$ induced by $e_{\pm}$. With respect to the splitting (\ref{ag}), the covariant tractor derivative $\nabla^{nc}$ on $\Gamma(\Lambda^{p+1,q+1}_{\mathcal{T}}(M))$ is given by\footnote{In comparison to \cite{nc} the roles of $\alpha_+$ and $\alpha_-$ are interchanged since the reference realises the parabolic subgroup $P$ as stabilizer of the line $\R e_+$.}
\begin{align*}
\nabla_X^{nc} \stackrel{g}{=} \begin{pmatrix}  \nabla^{g}_X &  - \iota_{X}  & X^{\flat} \wedge  & 0 \\ 
 -K^g(X)^{\flat} \wedge &  +  \nabla^{g}_X & 0 &  X^{\flat} \wedge  \\ 
 \iota_{K^g(X)} & 0 & \nabla^{g}_X &  \iota_X \\ 
 0 &  \iota_{K^g(X)}  & K^g(X)^{\flat} \wedge & \nabla^{g}_X  \end{pmatrix}.
\end{align*}
Using this expression, it is straightforward to calculate that if $\alpha$ is a normal twistor $(k+1)$-form, then $\alpha_+$ is a nc-Killing $k$-form and $\alpha_0 , \alpha_{\mp}, \alpha_-$ are uniquely determined by $\alpha_+$. 	On the other hand, given a nc-Killing $k-$form $\alpha \in \Omega^k(M)$ wrt. the metric $g$, there is a unique twistor $(k+1)-$form $\alpha^{k+1} \in \Omega^{k+1}_{\Pe^1}(M)$ such that $\left(\alpha^{k+1}\right)_+ = \alpha$ holds. Now let $\ph \in \Gamma(S^g)$ be a twistor spinor with associated parallel spin tractor $\psi \in \Gamma(\Se)$. Then it is straightforward to calculate that wrt. $g \in c$
\begin{align}
\left(\alpha^{k+1}_{\psi} \right)_+ = d_1 \cdot \alpha^k_{\ph}\text{ and }
\left(\alpha^{k+1}_{\psi} \right)_- = d_2 \cdot \alpha^k_{D^g\ph} \label{77}
\end{align}
where $d_{1,2}$ are nonzero constants not depending on $\psi$ (cf.\cite{ns}). Moreover, it holds that
\begin{align}
\ph(x) = 0 &\Rightarrow \alpha_{\psi}^{k+1}(x) = d_2 s_-^{\flat}(x) \wedge \alpha^k_{D^g \ph}(x), \label{ph}\\
D^g \ph(x) = 0 &\Rightarrow \alpha_{\psi}^{k+1}(x) = d_1 s_+^{\flat}(x) \wedge \alpha^k_{\ph}(x). \nonumber
\end{align}
Note that these formulas determine the $SO(p+1,q+1)$ orbit type of the parallel form $\alpha_{\psi}^{k+1}$ \textit{on all of M}. The existence of certain normal twistor forms has many interesting implications on the (local) geometry of $M$ as studied in \cite{leihabil}.
To summarize, one has the following implications between the objects introduced so far:
\begin{align*}
\begin{xy}
  \xymatrix{
      \ph \in \Gamma(S^g) \text{TS} \ar@{<->}[r]^{g \in c} \ar[dd]^{\text{nc-Killing}}    &   \psi \in \Gamma(\mathcal{S}_{\mathcal{T}}(M))\ar[r]^{\text{ Hol-Pr.}} \ar[dd]^{\text{normal twistor}} & v_{\psi} \in  \Delta_{p+1,q+1} \ar[dd]^{\text{alg. Dirac}} \ar[rd]^{\text{cond. for}}& \\
      & & & Hol(M,c) \\
      \alpha^p_{\ph} \in \Omega^p(M) \ar@{<->}[r]^{g \in c}             &  \alpha^{p+1}_{\psi} \in \Omega^{p+1}_{\Pe^1}(M)\ar[r]^{\text{ Hol-Pr.}}  & \alpha^{p+1}_{v_{\psi}} \in \Lambda^{p+1}_{p+1,q+1} \ar[ru]^{\text{cond.for}}& \\
  }
\end{xy}
\end{align*}

\section{Geometries admitting totally lightlike, holonomy-invariant subspaces} \label{s3}
Let $\psi \in \Gamma(\mathcal{S}(M))$ be a parallel spin tractor. We set $\mathcal{H}_{\psi(x)}:= \{ v \in \mathcal{T}_x(M) \mid  v \cdot \psi(x) = 0 \}$. This leads to a totally lightlike and parallel distribution $\mathcal{H}_{\psi} \subset \mathcal{T}(M)$. We want to prove that the twistor spinor induced by $\psi$ via fixing a metric in the conformal class is locally equivalent to a parallel spinor iff $\mathcal{H}_{\psi}$ is nontrivial. Main ingredient is the following statement about totally lightlike parallel distributions in the standard tractor bundle:
\begin{Proposition} \label{ct}
Let $(M,c)$ be a conformal manifold of dimension $n \geq 3$ and let $\mathcal{H} \subset \mathcal{T}(M)$ be a totally lightlike distribution of dimension $k \geq 1$ which is parallel wrt. the Cartan connection $\nabla^{nc}$. Then there is an open, dense subset $\widetilde{M} \subset M$ such that for every point $x \in \widetilde{M}$ there is an open neighborhood $U_x \subset \widetilde{M}$ and a metric $g \in c_{|U_x}$ such that wrt. the splitting (\ref{ide}) $\mathcal{H}$ is locally given by
\begin{align*}
\mathcal{H}_{|U_x} \stackrel{g}{=} \text{span }\left(\begin{pmatrix} 0 \\ K_1\\ 0 \end{pmatrix},...,\begin{pmatrix} 0 \\ K_{k-1}\\ 0 \end{pmatrix}, \begin{pmatrix} 0 \\ 0\\ 1 \end{pmatrix}\right)
\end{align*}
for lightlike vectorfields $K_i \in \mathfrak{X}(U_x)$.
\end{Proposition}

\begin{pf}If $k=1$, this is a well known fact (cf. \cite{baju}). We can adopt parts of the (first steps of the) proof and the notation from \cite{ln} where the statement is proved for $k=2$ and we may then also assume that $k>2$. However, we later use a different method. To start with, we set $\mathcal{L}:=\mathcal{I}^{\bot} \cap \mathcal{H}$, where $\mathcal{I}$ is the isotropic line defining the parabolic subgroup $P$. With respect to $g \in c$ one has that $\mathcal{L} = \left\{ X \in \mathcal{H} \mid X = \begin{pmatrix} \alpha \\Y \\ 0 \end{pmatrix} \right\}$. Note that $L:= \text{pr}_{TM} \mathcal{L} \subset TM$ is conformally invariant.\\

\textit{Step 1:}\\
We claim that there is an open, dense subset\footnote{In this proof, in order to keep notation short, whenever we restrict to an open, dense subset we again call it $M$.} $\widetilde{M} \subset M$ such that $\text{rk }\mathcal{L}_{|\widetilde{M}} = k-1$:
Note that $\mathcal{L} \neq \{0 \}$ as otherwise $\mathcal{H}$ would have rank 1. Moreover, there is no open set in $M$ on which $\text{rk }\mathcal{L} = k$ as follows from Lemma 2 in \cite{ln}. Consequently, there is an open, dense subset (which we again call $M$) over which $0 < \text{rk } \mathcal{L} < k$. Now let $x \in M$ and fix a basis $L_1,...,L_s$ of $\mathcal{L}_x$, where $s \leq k-1$. We may add tractors $Z_l = \begin{pmatrix} a_l \\ Y_l \\ 1 \end{pmatrix}$ for $1 \leq l \leq k-s$ such that $L_1,...,L_s,Z_1,...,Z_{k-s}$ is a basis of of $\mathcal{H}_x$. We know that $k-s \geq 1$. If $k-s >1$ we may form new basis vectors $Z_1+Z_2$ and $Z_1-Z_2$. However, $Z_1-Z_2 \in \mathcal{L}_x$. Thus, $k-s=1$, which shows that rk $\mathcal{L}_x = k-1$.\\
\newline
\textit{Step 2:}\\
We claim that also $L=\text{pr}_{TM}\mathcal{L}$ has rank $k-1$ locally around each point $x \in M$. To this end, let $g \in c$ be arbitrary. Then we choose generators of $\mathcal{L}$ around $x$ such that locally $\mathcal{L} \stackrel{g}{=} \text{span } \left( \begin{pmatrix} \sigma_1 \\ \widetilde{K}_1 \\ 0 \end{pmatrix},...,\begin{pmatrix} \sigma_{k-1} \\ \widetilde{K}_{k-1} \\ 0 \end{pmatrix}\right)$. We may assume that $\sigma_1(x) \neq 0$. Otherwise, we find $\ph \in C^{\infty}(M)$ with $\widetilde{K}_1(\ph)(x) \neq 0$ and consider the metric $\widetilde{g}=e^{2 \ph}g$ instead (cf. (\ref{ide})). Moreover, we may by rescaling the generators assume that there is a neighborhood $U$ of $x$ on which $\sigma_1 \equiv 1$ and $|\sigma_i| < 1$ for $i=2,...,k-1$. By linear algebra we then see that there are lightlike vectorfields $K_i$ for $i=1,...,k-1$ such that wrt. ${g}$ on $U$ 
\begin{align}
\mathcal{L} \stackrel{{g}}{=} \text{span } \left( \begin{pmatrix} 1 \\ {K}_1 \\ 0 \end{pmatrix}, \begin{pmatrix} 0 \\ {K}_2 \\ 0 \end{pmatrix},...,\begin{pmatrix} 0 \\ {K}_{k-1} \\ 0 \end{pmatrix}\right). \label{dd}
\end{align}
Suppose now that there is an open set on which  $\begin{pmatrix} 1 \\ 0 \\ 0 \end{pmatrix} \in \mathcal{L}$. Differentiating in direction $X \in TM$ yields that
$\nabla^{nc}_X \begin{pmatrix} 1 \\ 0 \\ 0 \end{pmatrix} = \begin{pmatrix} 0 \\ X \\ 0 \end{pmatrix} \in \mathcal{H}$ for all $X$ as $\mathcal{H}$ is parallel. This is not possible for dimensional reasons. Consequently, on an open and dense subset the vectors $K_1,...,K_{k-1}$ are linearly independent and as 
$L= \text{span} (K_1,...,K_{k-1})$ this shows that there is an open and dense subset of $M$ on which the rank of $L$ is maximal.\\
\newline
\textit{Step 3:}\\
It follows precisely as in the $k=2$-case from \cite{ln}, Lemma 3 and 4, that $L^{\bot} = \text{pr}_{TM} \left( \mathcal{H}^{\bot} \cap \mathcal{I}^{\bot} \right)$\\
\newline
\textit{Step 4:}\\
In the setting of Step 2 we consider the local representation (\ref{dd}) of $\mathcal{L}$ wrt. ${g}$ and set $L':=\text{span }(K_2,...,K_{k-1})$. Both $L$ and $L'$ are integrable distributions: Let $i,j \in \{2,...,k-1\}$. As $\mathcal{H}$ is parallel and totally lightlike we have that 
$\nabla_{K_i} \begin{pmatrix} 0 \\ K_j \\ 0 \end{pmatrix} = \begin{pmatrix} -P^g (K_i,K_j) \\ \nabla^g_{K_i}K_j \\ -g(K_i,K_j) \end{pmatrix} \in \Gamma(\mathcal{L})$. 
Switching the roles of $i$ and $j$ and taking the difference yields $\begin{pmatrix} 0 \\ \left[K_i,K_j\right] \\ 0 \end{pmatrix} \in \Gamma(\mathcal{L})$. Thus $\left[K_i,K_j\right] \in L'$. Similarly one shows that even
\begin{align}
[K_1,L'] \subset L' \label{ic}
\end{align}
\textit{Step 5:}\\
We now apply Frobenius Theorem: For every (fixed) point $y$ of (an open and dense subset of ) $M$ we find a local chart $(U,\ph=(x_1,...,x_n))$ centered at $y$ with $\ph(U) = \{(x_1,...,x_n) \in \R^n \mid |x_i| < \epsilon \}$ such that the leaves $A_{c_k,...,c_n}=\{a \in U \mid x_k(a) = c_1,...,x_n(a)=c_n \} \subset U$ are integral manifolds for $L$ for every choice of $c_j$ with $|c_j| <\epsilon$. It holds that $L_U = \text{span }\left(\frac{\partial}{\partial x_1},...,\frac{\partial}{\partial x_{k-1}}\right)$ and moreover the coordinates may be chosen such that $K_1 = \frac{\partial}{\partial x_1}$ over $U$. After applying some linear algebra to the generators of $L'$ and restricting $U$ if necessary, we may assume that generators of $L'$ are given on $U$ by
\begin{align}
K_{i \geq 2} =  \alpha_i \frac{\partial}{\partial x_1} + \frac{\partial}{\partial x_i} \label{ki}
\end{align}
for certain smooth functions $\alpha_i \in C^{\infty}(U)$. The integrability condition (\ref{ic}) implies that
\begin{align}
\frac{\partial}{\partial x_1} \alpha_i = 0 \text{ for }i=2,...,k-1. \label{tt}
\end{align}
The integrability of $L'$ and (\ref{tt}) then yield that
\begin{align}
\frac{\partial}{\partial x_i} \alpha_j - \frac{\partial}{\partial x_j} \alpha_i = 0 \text{ for }i,j=2,...,k-1. \label{jj}
\end{align}
For fixed $c_k,...,c_n$ as above we consider the submanifold $A_{c_k,...,c_n}$ and the differential form $\alpha := - \sum_{i=1}^{k-1} \alpha_i dx_i \in \Omega^1\left(A_{c_k,...,c_n} \right),$ where the $\alpha_{i \geq 2}$ are restrictions of the functions appearing in (\ref{ki}) and we set $\alpha_1 \equiv -1$. (\ref{tt}) and (\ref{jj}) yield that $d\alpha = 0$. Whence there exists by the Poincare Lemma (after restricting $U$ if necessary) a unique $\sigma_{c_k,...,c_n} \in C^{\infty}\left(A_{c_k,...,c_n} \right)$ with $\sigma_{c_k,...,c_n}(\ph^{-1}(0,...,0,c_k,...,c_n)) = 0$ and
\begin{align}
\frac{\partial}{\partial x_1} \sigma_{c_k,...,c_n} &= 1, \\
\frac{\partial}{\partial x_i} \sigma_{c_k,...,c_n} &= - \alpha_i \text{ for }i=2,...,k-1.
\end{align}
We then define $\sigma \in C^{\infty}(U)$ via $\sigma(\ph^{-1}(x_1,....,x_n)):=\sigma_{x_k,...,x_n}(\ph^{-1}(x_1,...,x_{n}))$ and observe that on $U$
\begin{align}
\frac{\partial}{\partial x_1} \sigma = 1\text{, }\frac{\partial}{\partial x_i} \sigma = - \alpha_i \text{ for }i=2,...,k-1. \label{hhh}
\end{align}
\textit{Step 6:}\\
(\ref{ki}) and (\ref{dd}) imply that on $U$ $K_1(\sigma)=1$ and $K_i(\sigma) = 0$ for $i=2,...,k-1$. We now consider the rescaled metric ${\widetilde{g}}=e^{2\sigma}{g}$ on $U$. The transformation formula (\ref{tra}) and (\ref{hhh}) then show that wrt. this metric $\mathcal{L}$ is given by
\begin{align*}
\mathcal{L}_U = \text{span} \left(\begin{pmatrix} 0 \\ K_1 \\ 0 \end{pmatrix},...,\begin{pmatrix} 0 \\ K_{k-1} \\ 0 \end{pmatrix} \right).
\end{align*}
We may add one generator $\begin{pmatrix} \beta \\ X \\1 \end{pmatrix} \in \Gamma(U,\mathcal{H})$ such that pointwise (wrt. ${\widetilde{g}}$) $\mathcal{H}= \mathcal{L} \oplus \text{span }\begin{pmatrix} \beta \\ X \\1 \end{pmatrix}$. It follows that $X \in L^{\bot}$. By step $3$ there exists a smooth function $b$ on $U$ with $X= \text{pr}_{TM} \begin{pmatrix} b \\ X \\0 \end{pmatrix}$ and $\begin{pmatrix} b \\ X \\0 \end{pmatrix} \in \mathcal{H}^{\bot}$. As $\mathcal{H}$ is lightlike, it follows that $b=\beta$. Therefore we have that $\begin{pmatrix} 0 \\ 0 \\1 \end{pmatrix} \in \mathcal{H}^{\bot}$ over $U$. However, this implies that $\beta = 0$. $\text{dim }\mathcal{H}+\text{dim }\mathcal{H}^{\bot}=n+2$ and dimension count show that $X$ has to be a linear combination of the $K_i$, $i=1,...,k-1$, and passing to new generators then proves the Proposition.
\end{pf}

We study some consequences. In the setting of Proposition \ref{ct} we have that $\mathcal{H}$ is parallel iff $\mathcal{H}^{\bot}$ is parallel. Locally, we have wrt. the metric appearing in Propostion \ref{ct} that $\mathcal{H}^{\bot}=\text{span} \left( \begin{pmatrix} 0 \\ X \\ \sigma \end{pmatrix} \mid X \in L^{\bot} \right)$. It follows that $\mathcal{H}^{\bot}$ is parallel iff

\begin{align*}
\nabla^{nc}_Y \begin{pmatrix} 0 \\ X \\ \sigma \end{pmatrix} = \begin{pmatrix} -P^g(X,Y) \\ \nabla^g_Y X + \sigma P^g(Y) \\ -g(X,Y) \end{pmatrix} \in \Gamma(U,\mathcal{H}^{\bot})
\end{align*}
for all $X \in \Gamma(U,L^{\bot})$ and $Y \in \mathfrak{X}(U)$. This is equivalent to parallelity of $L$,  $P^g(X)=0$ for all $X \in L^{\bot}$ and $P^g(TM) \subset L$. As in \cite{lst} we conclude that the scalar curvature is zero. Thus we have proved the following Proposition.

\begin{Proposition} \label{gg}
If on a conformal manifold $(M,c)$ there exists a totally lightlike, $k$-dimensional parallel distribution in $\mathcal{T}(M)$, then every point of some open and dense subset admits a neighborhood $U$, a metric $g \in c_U$ and a $k-1$-dimensional totally lightlike distribution $L \subset TU$ such that
\begin{align}
Ric^g(TU) \subset L,\text{ and }L \text{ is parallel wrt. }\nabla^g. \label{4}
\end{align}
Conversely, if $U \subset M$ is an open set equipped with a metric $g \in c_U$ and a $k-1$-dimensional totally lightlike distribution $L \subset TU$ such that (\ref{4}) holds, then $L$ gives rise to a $k-$dimensional totally lightlike, parallel distribution $ \begin{pmatrix} 0 \\ L \\ 0 \end{pmatrix} \oplus \text{span }\begin{pmatrix} 0 \\ 0 \\1 \end{pmatrix}$ in $\mathcal{T}(U)$.
\end{Proposition}
In case $k=1$, this means that there is locally a Ricci-flat metric in the conformal class. In case $k=2$ this describes conformal pure radiation metric with parallel rays as discussed in \cite{ln}. Proposition \ref{gg} also generalizes results from \cite{nc} where the statement is proved under the additional condition that the totally lightlike distribution arises from a decomposable, totally lightlike twistor $k-$form. One proves precisely as in \cite{ln}, Remark 2, that in the setting of Proposition \ref{ct} one gets the conformally invariant curvature condition
\[ W^g(L,L^{\bot},\cdot,\cdot) = 0 \]
for the Weyl tensor for the existence of a totally lightlike, parallel null-plane in the tractor bundle.\\
\newline
We apply these results to the case of twistor spinors on conformal spin manifolds. Let $\psi \in \Gamma(\mathcal{S})$ be a parallel spin tractor on $(M^{p,q},c)$ and for $g \in c$ let $\ph \in \Gamma(S^g)$ be the associated twistor spinor. As $\psi$ is parallel, the pointwise kernel of Clifford multiplication ker $\psi(x)=\{v \in \mathcal{T}_x(M) \mid v \cdot \psi(x)= 0\}$ yields a totally lightlike and parallel distribution $\mathcal{H}_{\psi} \subset \mathcal{T}(M)$. Similarly, if even $\ph$ is parallel, we get a totally lightlike, parallel distribution $L_{\ph} \subset TM$.
One then has the following immediate consequence from Proposition \ref{ct}:

\begin{Proposition} \label{tms}
If $\psi \in \Gamma(\mathcal{S})$ is a parallel spin tractor with $\mathcal{H}_{\psi} \neq 0$, then there is an open and dense subset $\widetilde{M} \subset M$ such that on $\widetilde{M}$ the associated twistor spinor is locally conformally equivalent to a parallel spinor.
\end{Proposition}
\begin{pf}
We notice that Proposition \ref{ct} yields the desired $\widetilde{M}$ and for $x \in \widetilde{M}$ a neighborhood $U$ and a local metric $g_{U}\in c_U$ such that $s_+ \in \mathcal{H}_{\psi_U}$. If we decompose $\psi$ on $U$ wrt. the metric $g$ as in Theorem \ref{90}, i.e. $\psi_{|U} = \left[\left[\sigma^g(l),e \right], e_-w + e_+w \right]$ for some function $w:U \rightarrow \Delta_{p+1,q+1}$, it follows that $e_+e_-w=0$ on $U$  which implies that $e_-w=0$. However, by Theorem \ref{90} it follows that on $U$ we have $D^g_{\ph} = \text{proj}_-^g (\psi) = 0$. Thus, $\ph$ is on $U$ both harmonic and a twistor spinor and therefore parallel wrt. $g$.
\end{pf}

Note that by the same argumentation as in the last proof every parallel spinor $\Gamma(S^g) \ni \ph \stackrel{g}{\leftrightarrow} \psi \in \Gamma(\mathcal{S})$ satisfies $\mathcal{H}_{\psi} \neq 0$. Whence, Proposition \ref{tms} yields locally an \textit{equivalent} characterisation of those parallel spin tractors which lead to parallel spinors in terms of conformally invariant objects. In terms of the original data the Proposition can be rephrased as follows: Note that wrt. the decompositon (\ref{gdg}) the requirement $\mathcal{H}_{\psi} \neq 0$ is equivalent to say that there is $x \in M$, $g \in c$ and a nontrivial tripel $(\alpha,X,\beta) \in \R \oplus T_xM \oplus \R$ such that 
\begin{align*}
X \cdot \ph(x) + \alpha D^g\ph(x)=0,\\
X \cdot D^g\ph(x) + \beta \ph(x) = 0.
\end{align*}
Thus, if, for example, $D^g\ph$ vanishes at \textit{some} point for some metric in the conformal class, then the twistor spinor is already locally equivalent to a parallel spinor locally around every point (up to a singular set). Moreover, Proposition \ref{tms} admits several further consequences and applications which contribute to the classification problem for local geometries admitting twistor spinors on pseudo-Riemannian manifolds.\\
We first describe how it is related to and generalizes other results obtained for the Riemannian and Lorentzian case. For a Riemannian spin manifold $(M^n,g)$ with twistor spinor $\ph$ one has that $\left(M \backslash Z_{\ph},\widetilde{g}=\frac{1}{||\ph||^4}\right)$ is an Einstein space of nonnegative scalar curvature $\widetilde{R}$. If $\widetilde{R}>0$, then the rescaled spinor decomposes into a sum of two Killing spinors whereas in case $\widetilde{R}=0$ one has a Ricci-flat metric with parallel spinor. Proposition \ref{tms} precisely describes the last case in which dim $\mathcal{H}_{\psi}=1$. For the Lorentzian case, Lemma \ref{0} yields a relation between Proposition \ref{tms} and the classificaion of twistor spinors on Lorentzian manifolds using the nc-Killingform theory in \cite{leihabil}.
\begin{satz}[\cite{leihabil}; Thm.10] \label{bg} 
Let $\varphi \in \Gamma(S^g_{\C})$ be a spinor on a Lorentzian spin manifold of dimension $n\geq 3$. Then one of the following holds on an open and dense subset $\widetilde{M} \subset M$:
\begin{enumerate}
\item $\alpha_{\psi}^2=l_1^{\flat} \wedge l_2^{\flat}$ for standard tractors $l_1,l_2$ which span a totally lightlike plane.\\
In this case, $\ph$ is locally conformally equivalent to a parallel spinor with lightlike Dirac current $V_{\ph}$ on a Brinkmann space.
\item $\alpha_{\psi}^2=l^{\flat} \wedge t^{\flat}$ where $l$ is a lightlike, $t$ is an orthogonal, timelike standard tractor.\\
$(M,g)$ is locally conformally equivalent to $(\mathbb{R},-dt^2) \times (N_1, h_1) \times \cdots \times (N_r, h_r)$, where the $(N_i,h_i)$ are Ricci-flat K\"ahler, hyper-K\"ahler, $G_2$-or $Spin(7)$-manifolds.
\item $\alpha_{\psi}^2$ is of Kaehler-type at every point (cf. Remark \ref{rem}).\\
The following cases can occur:
\begin{enumerate}
\item The dimension $n$ is odd and the space is locally equivalent to a Lorentzian Einstein-Sasaki manifold.
\item $n$ is even and $(M,g)$ is locally conformally equivalent to a Fefferman space.
\item $\alpha_{\psi}^2$ is a volume form on a nondgenerate subbundle $\mathcal{V} \subset \mathcal{T}(M)$. Then there exists locally a product metric $g_1 \times g_2 \in [g]$ on $M$, where $g_1$ is a Lorentzian Einstein-Sasaki metric on a space $M_1$ of dimension $n_1=2 \cdot rk(\alpha_1(\ph))+1$ admitting a Killing spinor and $g_2$ is a Riemannian Einstein metric with Killing spinor on a space $M_2$ of positive scalar curvature $scal^{g_2} = \frac{(n-n_1)(n-n_1-1)}{n_1 (n_1 - 1)}scal^{g_1}$.
\end{enumerate}
\end{enumerate} 
\end{satz}
Lemma \ref{0} shows us that $\mathcal{H}_{\psi} \neq 0$ occurs exactly in the first two cases in which we get a parallel spinor as also follows from Proposition \ref{tms}. In the third case, it holds that dim $H_{\psi}=0$ and the spinor cannot be rescaled to a parallel spinor. In particular, we have shown that the Killing spinors defining Lorentzian Einstein Sasaki structures (cf. \cite{boh}) can never be rescaled to parallel spinors.
\newline
We describe further geometric consequences implied by Propositon \ref{tms}. If $\ph$ can locally be rescaled to a parallel spinor, the vanishing of the torsion of $\nabla^g$ implies as a global consequence that $L_{\ph}$ is an \textit{integrable}, distribution on $\widetilde{M}$. Now fix $x \in \widetilde{M}$ and let $U \subset \widetilde{M}$ be an open neighborhood with metric $g \in c_{|U}$ such that $\ph$ is parallel wrt. $g$ on $U$. One has that  $Ric^g(TU) \subset L_{\ph|U}$ as implied by Proposition \ref{gg}. However, note that this also follows from the well-known fact that $Ric^g(TM) \cdot \ph = 0$ for any parallel spinor $\ph$. Moreover, it follows from Lemma \ref{0} or Proposition \ref{ct} that $\text{dim }H_{\psi}=\text{dim }L_{\ph_{|U}}+1$. In case that $k:=\text{dim }L_{\ph_{|U}}>0$, $Hol(U,g)$ acts reducible with a fixed totally lightlike $k-$dimensional subspace. If $k=p$, i.e. $\ph$ is a pure spinor on $U$, it follows from Lemma \ref{0} that even a totally isotropic $p-$form is fixed. If $k=p-1$, $Hol(U,g)$ fixes a $p-$form of type $\alpha^p_{\ph}=l_1^{\flat} \wedge...\wedge l_{p-1}^{\flat} \wedge t^{\flat}$, where $t$ is not lightlike and it follows that even then totally lightlike form  $l_1^{\flat} \wedge...\wedge l_{p-1}^{\flat}$is fxed by the holonomy representation. If $k=0$ it follows from Proposition \ref{tms} or $Ric^g(X)\cdot \ph=0$ that $g$ is a Ricci-flat metric on $U$. There is a complete list of possible irreducible, non locally symmetric holonomy groups for this case as to be found in \cite{kath}.
We summarize these results as follows:

\begin{Proposition} \label{trust}
Let $\psi$ be a parallel spin tractor with $\mathcal{H}_{\psi} \neq 0$. Then there is an open, dense subset $\widetilde{M} \subset M$ such that $L_{\ph}$ is an integrable distribution on $\widetilde{M}$. Moreover, any $x \in \widetilde{M}$ admits an open neighborhood $U \subset \widetilde{M}$ and a metric $g \in c_{|U}$ such that $\ph$ is a parallel spinor on $(U,g)$ and one of the following cases occurs:
\begin{enumerate}
\item $k:=\text{dim }L_{\ph} \neq 0$. In this case, $Hol(U,g)$ acts reducible with fixed $k-$dimensional totally lightlike subspace $L$ and $Ric^g(TU)\subset L$. Moreover, if $k=p,p-1$ there is a totally isotropic parallel $k-$form. 
\item $k:=\text{dim }L_{\ph} = 0$. The space $(U,g)$ is Ricci-flat. If it is not locally symmetric and $Hol(U,g)$ acts irreducible, then it is one of the list in \cite{kath}.
\end{enumerate}
\end{Proposition}
We further remark that similar integrability conditions for \textit{pure} twistor spinors have been derived in \cite{tag1} and \cite{tag2}. In split signature $(m+1,m)$ where $\Delta_{m+1,m}^{\C}$ admits a real structure and one can speak about real spinor fields one can say even more about parallel pure spinor fields by using results from \cite{kath} which give an explicit normal form for the metric for this case. More concretely, let $(M,h)$ be a pseudo-Riemannian spin manifold of split signature $(m+1,m)$ admitting a real pure parallel spinor field. Then one can find for every point in $M$ local coordinates  $(x,y,z)$ , $x=(x_1,...,x_m)$, $y=(y^1,...,y^m)$ around this point such that 
\begin{align} h = -dz^2 - 4 \sum_{i=1}^m dx_i dy^i - 4 \sum_{i,j=1}^m g_{ij} dy^i dy^j, \label{purre1} \end{align}
where $g_{ij}$ are functions depending on $x,y$ and $z$ and satisfying
\begin{align}g_{ij} = g_{ji} \text{ for }i,j=1,...,m \text{,  }\sum_{i=1}^m \frac{\partial g_{ik}}{\partial x_i} = 0 \text{ for } k=1,...,m. \label{pure2} \end{align}
Conversely, if one uses (\ref{purre1}) and (\ref{pure2}) to define a metric $h$ on a connected open set $U \subset \R^{2m+1}$, then $(U,h)$ is spin and admits a real pure parallel spinor. Similar statements hold in case $(p,q)=(m,m)$, where one has to omit the last coordinate etc.\\

As a special application we consider twistor spinors equivalent to parallel spinors in case $p=2$. If $\mathcal{H}_{\psi} \neq \{0\}$, then in the above notation one has a parallel 2-form $\alpha_{\ph}^2$ on $(U,g)$. The $SO^+(2,n-2)$-orbit type of this form must be one of the list from Remark \ref{rem}: The first form, $\alpha_{\ph}^2 = l^{\flat}_1 \wedge l^{\flat}_2$ corresponds to a parallel pure spinor. In the second case, $\alpha_{\ph}^2 = l^{\flat} \wedge t^{\flat}$, we can conlcude that there is a nontrivial lightlike, parallel vectorfield and thus $(U,g)$ is a Brinkmann space. In the third case, $(U,g)$ is Ricci-flat (as $L_{\ph} = 0$) and $Hol(U,g)$ leaves invariant a (possibly trivial) $n-2m$ dimensional nondgenerate subspace $E^{\bot}$ and $\alpha_{\ph}^2$ is Kaehler on $E$. It follows with Remark (\ref{rem}) that there is a local splitting $(U,g) \cong (U_1,g_1) \times (U_2,g_2)$, where the first factor is Ricci-flat pseudo Kaehler of signature $(2,2m-2)$ and the second factor (which might be trivial) is Riemannian Ricci flat. Moreover, by Leitners argument from \cite{leihabil} both factors admit parallel spinors.

\section{The zero set of a twistor spinor}
In this section we want to describe the possible local shapes of the zero set $Z_{\ph}$ of a twistor spinor $\ph \in \Gamma(S^g)$ and study the properties and related local geometries of the spinor off its zero set. It is shown in \cite{bfkg} that in the Riemannian case the zero set consists of a countable union of isolated points. For the Lorentzian case, \cite{ns} shows that $Z_{\ph}$ - if nonempty - consists either of isolated points or of isolated images of lightlike geodesiscs. Moreover, one has that for a given $x \in Z_{\ph}$, there is an open neighborhood $U$ of $x$ in $M$ and $V$ of $0$ in $T_xM$ such that
\begin{align}
Z_{\ph} \cap U = \text{exp}_x (\text{ker }D^g \ph(x) \cap V). \label{ffzz}
\end{align}
The proof of (\ref{ffzz}) relies on the investigation of the zero set of $V_{\ph}$, being a conformal vector field which additionally satisfies $\iota_{V_{\ph}}W^g = 0$. We show that (\ref{ffzz}) holds in all signatures by making use of the holonomy reduction procedure for general Cartan geometries as introduced in \cite{cgh}. Applied to our setting, this reads as follows: Let $\psi \in \Gamma(\mathcal{S}(M))$ be a $\nabla^{nc}$-parallel spin tractor. We view $\mathcal{S}(M) = \overline{\mathcal{Q}^1_+}\times_{\widetilde{G}}\Delta_{p+1,q+1}$. By standard principle bundle theory, $\psi$ then corresponds to a $\widetilde{G}-$equivariant smooth map $\widehat{\psi}: \overline{\mathcal{Q}^1_+} \rightarrow \Delta_{p+1,q+1}$. As $\psi$ is parallel, the image $\mathcal{O}:=\widehat{\psi}\left( \overline{\mathcal{Q}^1_+} \right) \subset \Delta_{p+1,q+1}$ is a orbit wrt. the $\widetilde{G}$-action, called the $\widetilde{G}$\textit{-type} of $\psi$. We now bring into play that $\nabla^{nc}$ is induced by $\left(\mathcal{Q}^1_+,\widetilde{\omega}^{nc}\right)$: Let $x \in M$. We define the $\widetilde{P}-$\textit{type} of $x$ wrt. $\psi$ to be the $\widetilde{P}-$orbit $\widehat{\psi} \left(\mathcal{Q}^1_+ \right) \subset \mathcal{O} \subset \Delta_{p+1,q+1}$ which may change over $x \in M$. $M$ then decmposes into a disjoint union according to $\widetilde{P}-$types, each of which is an initial submanifold of $M$. Then Proposition 2.7 from \cite{cgh} applied to our setting immediatly yields the following:

\begin{Proposition}\label{ddg}
Let $(M^{p,q},c)$ be a conformal spin manifold and let $\psi$ be a parallel spin tractor on $\left(\mathcal{Q}^1_+ \rightarrow M, \widetilde{\omega}^{nc} \right)$. For given $g \in c$ denote by $\ph \in \Gamma(S^g)$ the corresponding twistor spinor. Let $x \in M$. Then there is a parallel spin tractor $\phi$ on the homogeneous model $\left(\widetilde{G} \rightarrow \widetilde{G}/\widetilde{P} = \widehat{Q}^{p,q},\omega^{MC} \right)$ for which $x':=e\widetilde{P} \in \widetilde{G}/\widetilde{P}$ has the same $\widetilde{P}-$type wrt. $\phi$ that $x$ has wrt. $\psi$. Further, let $\ph'$ correspond to $\phi$ via a conformally flat metric ${g}_{St}$ on $\widehat{Q}^{p,q}$. Then there are open neighborhoods $N$ of $x$ in $M$ and $N'$ of $x'$ in $\widehat{Q}^{p,q}$ and a diffeomorphism $\Phi : N \rightarrow N'$ such that $\Phi(x) = x'$ and 
\begin{align*}
\Phi \left(Z_{\ph} \cap N \right) = Z_{\ph'} \cap N'.
\end{align*}
\end{Proposition}
As locally all possible shapes of the zero set already show up in the homogeneous model, we are led to study the zeroes of twistor spinors on $\widehat{C}^{p,q}$. Using (\ref{dat}) we identify $\widehat{Q}^{p,q}$ with the product $S^p \times S^q \subset \R^{p+1,q+1}$ equipped with the conformally flat metric $g_{St}:=-g_{S^p} + g_{S^q}$. We follow \cite{lei} in order to construct all twistor spinors on $\widehat{C}^{p,q}$. We decompose every $x \in \R^{n+2} \cong \R^{p+1} \times \R^{q+1}$ into $x=(x_1,x_2)$. There is a natural, globally defined orthonormal frame field on the normal bundle $N\widehat{Q}^{p,q}$, given by 
$\zeta_0(x)=(x_1,0)$ and $\zeta_{n+1}(x)=(0,x_2)$ for $x \in \widehat{Q}^{p,q}$. The spin structure on $\widehat{C}^{p,q}$ is then naturally induced by a standard spin structure on $\R^{p+1,q+1}$, and the spinor bundles are related by
\begin{align*}
S^{\R^{p+1,q+1}}_{|\widehat{Q}^{p,q}} \cong \underbrace{{Ann} \left(\zeta_0 + \zeta_{n+1} \right)}_{\cong S^{\widehat{Q}^{p,q},g}} \oplus \underbrace{{Ann}  \left(\zeta_0 - \zeta_{n+1} \right)}_{\cong S^{\widehat{Q}^{p,q},g}}.
\end{align*}
Wrt. this splitting, every spinor $\ph$ on $\R^{p+1,q+1}$ decomposes into $\ph=\ph_1+\ph_2$. For given $v \in \Delta_{p+1,q+1} \backslash \{0 \}$ we let $\ph_v(x):=x \cdot v$ for $x \in \R^{p+1,q+1}$, yielding a twistor spinor on $\R^{p+1,q+1}$. Using the relation between the spinor derivatives $\nabla^{\R^{p+1,q+1}}$ and $\nabla^{\widehat{Q}^{p,q}}$ it is straightforward to calculate that the induced spinor $\ph_{v,1}$ is a twistor spinor on $\left(\widehat{Q}^{p,q},g\right)$ (with $\ph_{v,2} \equiv 0$), and for dimensional reasons, all twistor spinors on the homogeneous model arise this way.
The next statement generalizes a classical result from \cite{li89} for the Riemannian case: A twistor spinor on the standard sphere admits at most one zero.
\begin{Proposition} \label{ho}
Let $\ph:=\ph_{v,1}$ be a twistor spinor on $(\widehat{Q}^{p,q},g=g_{St})$, induced by a twistor spinor $\ph_v$ on $\R^{p+1,q+1}$ as explained above. Suppose that there is $x \in Z_{\ph}$. Then it holds that
\begin{align*}
Z_{\ph} = \text{exp}_x \left(\text{ker }D^g \ph(x)\right)\text{ or }Z_{\ph}=\{x,-x \}.
\end{align*}
\end{Proposition}

\begin{pf}
First, one shows using the formulas in \cite{lei} that $D^g \ph_{v,1}(y)=n(-v+\frac{1}{2} \zeta_0 \cdot y \cdot v)$ for all $y \in \widehat{Q}^{p,q}$. In particular, $x \in Z_{\ph_{v,1}}$ implies that $\text{ker }D^g \ph (x) = \{t \in T_x \widehat{Q}^{p,q} \mid t \cdot v = 0 \}$. Now let $b \in T_x \widehat{Q}^{p,q} \backslash \{0 \}$ with $g_x(b,b)=\langle b,b \rangle_{p+1,q+1} =0$. One checks that the geodesic through $x$ in direction $b$ is given by $\delta_b(t)= \text{cos}(t ||b_1||)\cdot x +  \text{sin}(t||b_1||) \cdot \frac{b}{||b_1||}$ with $|| \cdot ||$ being the standard Euclidean norm on $\R^{p+1}$. If now additionally $b \cdot v = 0$, we have that $\delta_b(1) \cdot v = 0$ as $x \in Z_{\ph}$, i.e. $x \cdot v = 0$. This shows the $\supset$ direction. On the other hand, suppose that $y \in Z_{\ph}$. As $y \cdot v = x \cdot v = 0$, it follows that $0=(yx+xy)\cdot v = -2 \langle x,y \rangle_{p+1,q+1} v$. Since $\langle y_i,y_i \rangle = 1$ for $i=1,2$, we find $\alpha_i \in [0;\pi]$ and $d_1 \in \R^{p+1}, d_2 \in \R^{q+1}$ with $\langle x_i,d_i \rangle = 0$ and $||d_1||=||d_2||=1$ such that $y_i = \text{cos}(\alpha_i) \cdot x_i  + \text{sin}(\alpha_i) \cdot d_i $ for $i=1,2$. The condition $\langle x_1,y_1 \rangle = \langle x_2,y_2 \rangle $ then leads to $\alpha_1 = \alpha_2 = \alpha$. Thus, 
\begin{align*}
y = \text{cos} (\alpha) \cdot x  + \text{sin} (\alpha) \cdot d
\end{align*}
for $d=d_1+d_2 \in T_x\widehat{Q}^{p,q}$. If $\text{sin}(\alpha) \neq 0$, we conclude that $d \cdot v = 0$, and thus $d \in \text{ker }D^g \ph(x)$. As moreover $||d_1||=1$, we see that there is $t \in \R$ with $y=x \cdot \text{cos}(t||d_1||) + \frac{d}{||d_1||} \cdot \text{sin} (t ||d_1||) = \delta_d(t)=\delta_{td}(1)$.
If \text{sin}$(\alpha) = 0$, we have either that $y=x$ where the statement is trivial or $y=-x$. If $\text{ker }D^g \ph(x)$ is nontrivial in this situation, we may choose arbitrary $d \in \text{ker }D^g \ph(x) \backslash \{0 \}$ for a geodesic $\delta_{d}$ joining $x$ and $-x$. Otherwise $\text{ker }D^g \ph(x)=0$ and the situation $Z_{\ph} = \{x,-x \}$ occurs.
\end{pf}

The proof further yields the following for the flat model: Let $x' \in Z_{\ph'}$ and suppose that for some $w \in T_{x'}\widehat{Q}^{p,q} \cap W$, where $W:= \{ w \in T_{x'}\widehat{Q}^{p,q} \mid \sqrt{\langle w,w \rangle_{n+2}} < \frac{\pi}{2} \}$ it holds that $y=\text{exp}_{x'}(w) =\delta_{w}(1)\in Z_{\ph'}$. As $\langle y,y \rangle =0$ it follows that $w_1 \neq 0$ and $w_2 \neq 0$, and the geodesic $\delta_w$ is thus given by $\delta_w(t)=\text{cos}(t ||w_1||)\cdot x'_1 +  \text{sin}(t||w_1||) \cdot \frac{w_1}{||w_1||}+\text{cos}(t ||w_2||)\cdot x'_2 +  \text{sin}(t||w_2||) \cdot \frac{w_2}{||w_2||}$. Now $x',y \in Z_{\ph'}$ implies that $\langle x',\delta_{w}(1) \rangle_{p+1,q+1}=0$ which yields that $\text{cos}^2(||w_1||)=\text{cos}^2(||w_2||)$. However, $w \in W$ implies that $||w_1||=||w_2||$. Consequently, $\langle w,w \rangle_{p+1,q+1}=0$. Now $y \cdot v=x\cdot v=0$ leads to $w \cdot v =0$ as in the proof of the previous Proposition. This shows that $w \in \text{ker }D^g \ph'(x')$ . In the notation of Proposition \ref{ddg} we can therefore choose $N'=\text{exp}_{x'}(V)$ to be a sufficiently small normal neighborhood of $x'$ for some open neighborhood $V$ of $0$ in $T_{x'}\widehat{Q}^{p,q}$ with $V \subset W$, and get that 
\begin{align}
\Phi(Z_{\ph} \cap N) = Z_{\ph'} \cap N'=\text{exp}_{x'}\left(\text{ker } D^{g_{St}} \ph'(x') \cap V \right). \label{wich}
\end{align}
We now return to general twistor spinors on $(M^{p,q},c)$. We claim that in the notation of Proposition \ref{ddg} for the zero $x \in Z_{\ph}$ it holds that
\begin{align}
\text{dim ker }D^g \ph(x) = \text{dim ker }D^{g_{St}}(x'). \label{dasw}
\end{align}
Indeed, in the notation of Proposition \ref{ddg} and subsection \ref{tes} we have that 
\begin{align*}\psi(x)=\left[\sigma^g(l),e_- w \right] &\Rightarrow D^g \ph(x) = [l,\alpha(e_-w)] \\
\phi(x')=\left[\sigma^{g_{St}}(l'),e_- w' \right] &\Rightarrow D^{g_{St}} \ph'(x') = [l',\alpha(e_- w')]
\end{align*}
for spinors $w,w' \in \Delta_{p+1,q+1}$. As the $\widetilde{P}-$types coincide, there is $\widetilde{p} \in \widetilde{P}$ such that \begin{align}
\widetilde{p} \cdot (e_-w) = e_- \cdot w'. \label{pt}
\end{align}We therefore invetigate the $\widetilde{P}$-action on $Ann(e_-) \subset \Delta_{p+1,q+1}$ more closely. Consider the 2-fold covering $\lambda: Spin(p+1,q+1) \rightarrow SO(p+1,q+1)$ which is explicitly given by $\lambda(u)(x)=u \cdot x \cdot u^{-1}$ (cf. \cite{ba81}), i.e. \begin{align}\widetilde{p} \cdot x = \lambda(\widetilde{p})(x) \cdot \widetilde{p}. \label{ttt} \end{align}Now we can find $a \in \R^+$, $A \in O(p,q)$ and $v \in \left(\R \right)^*$ such that wrt. the splitting $\R^{p+1,q+1} \cong \R e_- \oplus \R^{p,q} \oplus \R e_+$ we have that $\lambda(\widetilde{p}) = \begin{pmatrix} a^{-1} & v & -\frac{1}{2} a \langle v,v \rangle_{p,q} \\ 0 & A & -aAv^{\sharp} \\ 0  & 0 & a \end{pmatrix}$. It is then a straightforward calculation using (\ref{ttt}) and the formulas for the action of $x \in \R^{n+1}$ on $Ann(e_-) \oplus Ann(e_+)$ that (wrt. appropriate bases) $\widetilde{p}$ acts as $\begin{pmatrix} X & Y \\ 0 & aX \end{pmatrix}$ on $Ann(e_-) \oplus Ann(e_+)$ for some $X \in GL(\Delta_{p,q})$ with \begin{align}\mu(Ax) \cdot X = X \cdot \mu(x) \label{X} \end{align} for $x \in \R^{p,q}$, where we identify $Ann(e_{\pm}) \cong \De$ as explained in (\ref{deco}).
(\ref{pt}) and (\ref{X}) then imply that $A( \text{ker }e_- w) = \text{ker }e_- w'$ which proves (\ref{dasw}).\\
Another interesting observation is that the quantity ker $D^g \ph(x)$ does not depend on the zero $x \in Z_{\ph}$ \footnote{Moreover, it does not depend on the chosen metric in the conformal class as can be seen directly from the transformation formulas.}. One way to see this is the structure of the parallel tractor form $\alpha_{\psi}^{p+1}$. We have already observed that 
\begin{align*}
 \alpha_{\psi}^{p+1}(x) = d_2 \cdot s_+^{\flat}(x) \wedge \alpha^p_{D^g \ph}(x) \text{ for } x \in Z_{\ph}.
 \end{align*}
Lemma \ref{0} then yields that dim ker $D^g \ph(x)= \text{dim ker }\psi(x)-1$ and the right side of this equation does not depend on $x \in Z_{\ph}$ as $\psi$ is parallel. 
The zero set $Z_{\ph}$ now turns out to be an embedded submanifold of $M$: Let $x \in Z_{\ph}$ be arbitrary. In the setting of Proposition \ref{ddg} and (\ref{wich}) we choose neighborhoods $N$ and $N'$ where we may assume that $N'=\text{exp}_{x'}(V)$ is a normal neighborhood of $x'$ as in (\ref{wich}). We then consider $\widetilde{\Phi}:= \left(\text{exp}_{x'}\right)_{|V}^{-1} \circ \Phi : N \rightarrow V$. Propositions \ref{ddg} and (\ref{wich}) yield that $\widetilde{\Phi} (Z_{\ph} \cap N) = \text{ker }D^g \ph'(x') \cap V$. We may compose this map with a linear isomorphism $A_{x'}:T_{x'}\widehat{Q}^{p,q} \rightarrow \R^n$ satisfying $A_{x'}(\text{ker }D^g \ph'(x')) = \R^k \times \{0 \}$, and in this way we obtain a submanifold chart for $Z_{\ph}$ around $x$.  This submanifold is totally lightlike, since for every curve $\gamma$ with $\text{Im }\gamma \subset Z_{\ph}$ the twistor equation yields that  $\gamma'(t) \cdot D^g \ph (\gamma(t)) = 0$. As $D^g \ph (\gamma(t)) \neq 0$, we have that $\gamma$ is isotropic.

In addition, Lemma 3.4.1 from \cite{lei} says in our notation that for every $x \in Z_{\ph}$ one has that $\text{exp}_x \left(\text{ker} D^g \ph(x) \cap D_x \right) \subset Z_{\ph}$, where $D_x$ is the maximal domain of definition for the exponential map at $x$. For dimensional reasons, one then has that $\text{exp}_x \left(\text{ker} D^g \ph(x) \cap V \right)$ is an open submanifold of the embedded submanifold $Z_{\ph} \cap U'$ for appropriate neighborhoods $V$ of $0 \in T_xM$ and $U'$ of $x \in M$. This yields (\ref{ffzz}) for arbitrary dimensions. We summarize our observations:

\begin{satz} \label{zss}
Let $\ph \in \Gamma(S^g)$ be a twistor spinor with zero. Then the zero set $Z_{\ph}$ is an embedded totally lightlike, totally geodesic submanifold of dimension $\text{ker }D^g \ph(x)$, where the last quantity does not depend on the choice of $x \in Z_{\ph}$. Moreover, for every $x \in Z_{\ph}$ there are open neighborhoods $U$ of $x$ in $M$ and $V$ of $0$ in $T_xM$ such that
\begin{align}
Z_{\ph} \cap U = \text{exp}_x \left(\text{ker }D^g \ph(x) \cap V \right). \label{ff}
\end{align}
\end{satz}
More loosely speaking, the connected components of the zero set consist either of an isolated point or of the image of a null-geodesic or of a totally null-plane etc. A mixture of two of these geometric objects cannot occur as the zero set of one single twistor spinor. The whole local geometry of the zero set is encoded in the quantity ker $D^g \ph(x)$. In case of a Ricci-parallel metric in the conformal class one has stronger results about the shape of the set $V$ appearing in (\ref{ff}) as explained in \cite{lei}.\\
\newline
We next show that the conformal class naturally induces a projective structure on the zero set of a twistor spinor. Recall that two connections $\nabla$ and $\widehat{\nabla}$ on a manifold $N$ are called projectively equivalent iff there exists a 1-form  $\Upsilon \in \Omega^1(N)$ such that 
\begin{align*} \widehat{\nabla}_X Y = \nabla_X Y + \Upsilon(Y) X + \Upsilon(X) Y \text{  } \forall X,Y \in \mathfrak{X}(N).\end{align*}
A more geometric interpretation is that two linear connections with the same torsion are projectively equivalent if and only if they admit the same geodesics as unparametrized curves. A projective structure on $N$ is an equivalence class of connections.

\begin{Proposition} \label{lps}
Let $\ph \in \Gamma(S^g)$ be a twistor spinor with $Z_{\ph} \neq \emptyset$ on $(M,c)$. Then for every $g \in c$ the Levi Civita connection $\nabla^g$ descends to a torsion-free linear connection $\nabla$ on $Z_{\ph}$. If $g$ and $\widetilde{g}$ are conformally equivalent, the induced connections $\nabla$ and $\widetilde{\nabla}$ are projectively equivalent, i.e., there is a natural construction
\begin{align*}
\ph \text{ on }(M,c) \rightarrow (Z_{\ph},\left[\nabla \right])
\end{align*}
from conformal structures and a twistor spinor with zero to torsion-free projective structures on the zero set.
\end{Proposition}

\begin{pf}
It follows directly from (\ref{ff}) that for $x \in Z_{\ph}$ the tangent space to the zero set is given by $T_x Z_{\ph}=\text{ker }D^g \ph(x) \subset T_xM$. In particular, $\text{ker }D^g \ph(x)$ does not depend on the choice of $g \in c$. For given $X,Y \in \mathfrak{X}(Z_{\ph})$ we then set $\nabla_X Y (x):= \left(\frac{\nabla^g}{dt}{\left(Y \circ \gamma \right)}\right)(0)$, where $\gamma :(-\epsilon, \epsilon) \rightarrow Z_{\ph} \subset M$ is the maximal geodesic in $M$ with $\gamma(0)=x$ and $\gamma'(0)=X(x)$. $\frac{\nabla^g}{dt}$ is the induced derivative along $\gamma$. We have to check that $\nabla_X Y \in \mathfrak{X}(Z_{\ph})$. As $(Y \cdot D^g \ph) \circ \gamma = 0$, it follows that
\begin{align*}
0 = \frac{\nabla^{S^g}}{dt} \left((Y \cdot D^g \ph) \circ \gamma \right)= \left(\left(\nabla_X Y \right) \cdot D^g \ph \right) \circ \gamma + (Y \circ \gamma) \cdot \underbrace{\frac{\nabla^{S^g}}{dt} (D^g \ph \circ \gamma )}_{=\frac{n}{2}K^g(X \circ \gamma) \cdot \left(\ph \circ \gamma \right) = 0}
\end{align*}
Consequently, $\nabla_X Y (x) \in \text{ker }D^g \ph(x) = T_xZ_{\ph}$. Clearly, this holds for every metric in the conformal class. The fact that $\nabla$ is torsion-free follows directly from the corresponding property of $\nabla^g$.
Now let $\widetilde{g}=e^{2 \sigma}g$ be a conformally equivalent metric. There is the well-known transformation formula
\begin{align*}
\nabla_X^{\widetilde{g}}Y = \nabla_X^g Y + X(\sigma)Y + Y(\sigma) X -g(X,Y) \text{grad}^g \sigma.
\end{align*}
As for $x \in Z_{\ph}$ the space $\text{ker }D^g \ph(x)$ is totally lightlike, it is a direct consequence of the definition of $\nabla$ that for all $X,Y \in \mathfrak{X}(Z_{\ph})$ we have 
\begin{align*}
\widetilde{\nabla}_X Y = \nabla_X Y + d\widehat{\sigma}(X) \cdot Y + d\widehat{\sigma}(Y) \cdot X,
\end{align*}
where $\widehat{\sigma}:=\sigma_{|Z_{\ph}}$. It follows that $\nabla$ and $\widehat{\nabla}$ are projectively equivalent.
\end{pf}

Note that as a direct consequence of the definitions it holds that $i^*R^g = R^{\nabla}$, where $i:Z_{\ph} \hookrightarrow M$ and $R^{\nabla}$ is the curvature tensor of the connection $\nabla$. In particular, if $c$ admits a flat representative then so does $[\nabla]$.\\
\newline
It is now natural to ask what can be said about the spinor and associated local geometries off the zero set if one knows the (local) structure of $Z_{\ph}$. In the Riemannian case, a twistor spinor is always parallel on a Ricci-flat space off the zero set.
For Lorentzian signature F. Leitner showed that in case of an isolated zero the Lorentzian metric is locally off the zero set isometric to a static monopole $-dt^2 + h$ where $h$ is a Riemannian Ricci-flat metric with parallel spinor. If the zero is not isolated, then off the zero set the space is locally conformally equivalent to a Brinkmann space with parallel spinor. Our results from section 3 show that in every signature the spinor is locally equivalent to a parallel spinor off the zero set. In fact, let $\psi \in \Gamma(\mathcal{S}(M))$ be a parallel spin tractor with associated twistor spinor $\ph \in \Gamma(S^g)$ for $g \in c$. Let $x \in Z_{\ph}$. It then holds \textit{at }$x$ that $\psi(x)=\left[\left[\sigma^g(l),e \right],0+e_-w\right]$ for some $w \in \Delta_{p+1,q+1}$. However, this means that $s_-(x) \in \text{ker }\psi(x)$. In particular, since the dimension of this kernel is constant over $M$, Proposition \ref{tms} applies and yields the next statement.

\begin{satz}
Let $\ph \in \Gamma(M,S^g)$ be a twistor spinor admitting a zero. Then there is an open dense subset $\widetilde{M} \subset M$ with $Z_{\ph} \subset M \backslash \widetilde{M}$ such that for every $x \in \widetilde{M}$ there is an open neighborhood $U_x \subset \widetilde{M}$ such that $\ph$ can be rescaled to a parallel spinor on $U_x$. 
\end{satz}

Our discussion from section \ref{s3} implies further consequences relating the shape of the zero set to local geometric structures off the zero set:

\begin{Proposition} \label{gogo}
Let $\ph \in \Gamma(S^g)$ be a twistor spinor with nonempty zero set $Z_{\ph}$. Then there is a set of singular points $sing(\ph) \subset M$ with $Z_{\ph} \subset sing(\ph)$ such that the following holds:
There is $0 \leq k \leq p$ such that $Z_{\ph}$ is an embedded $k-$dimensional totally lightlike submanifold. On $M \backslash sing(\ph)$, the spinor is locally conformally equivalent to a parallel spinor and the corresponding metric holonomy representation fixes a lotally lightlike subspace of dimension $k$. If $k=p$ or $k=p-1$ there is even a fixed totally lightlike $k-$form. If $k=0,$ i.e. the zero is isolated, there is locally a Ricci-flat metric in the conformal class.
\end{Proposition}
For the proof we observe first that for the number $k$ appearing in the Proposition it holds that $k= \text{dim ker }D^g \ph(x)$, where $x \in Z_{\ph}$. It follows in the notation of Lemma \ref{0} and Proposition \ref{trust} that $k= \text{dim }\mathcal{H}_{\psi}-1=\text{dim }L_{\ph|U}-1$ and $L_{\ph}$ is parallel and totally lightlike (wrt. to a suitable metric in the conformal class).\\
\newline
In case $p=2$, the discussion from the end of section \ref{s3} together with the last statement directly leads to the following relation between the shape of the zero set and local geometries:

\begin{Proposition}
Let $\ph \in \Gamma(S^g)$ be a twistor spinor with zero on $(M^{2,n-2},g)$. Then exactly one of the following cases occurs:
\begin{enumerate}
\item $Z_{\ph}$ consists locally of totally lightlike planes. In this case, the spinor is locally equivalent to a parallel spinor off the zero set and gives rise to a parallel totally lightlike 2-form.
\item $Z_{\ph}$ consists of isolated images of lightlike geodesics. In this case, the spinor is off the zero set locally conformally equivalent to a parallel spinor on a Brinkmann space.
\item $Z_{\ph}$ consists of isolated points. In this case there is for each point off the zero set an open neighborhood and a local metric in the conformal class such that the resulting space is isometric to a product $(U_1,g_1) \times (U_2,g_2)$ where the first factor is Ricci-flat pseudo-Kaehler and the second factor (which might be trivial) is Riemannian Ricci-flat. Both factors admit a parallel spinor.
\end{enumerate}
\end{Proposition}

\section{Low dimensions}

\subsection{Non-generic twistor spinors}
One important application of the main statement, Proposition \ref{tms}, is the case that $\psi$ is a pure spinor, i.e. dim $\mathcal{H}_{\psi}=p+1$. It immediatly follows with Lemma \ref{0} that $\alpha_{\psi}^{p+1}$ is totally lightlike in this case, and thus Proposition \ref{tms} and its Corollaries apply (if $p \neq0)$ yielding important consequences in small split-signatures due to the following mainly algebraic observations concerning the orbit structure of $\Delta_{p,q}$ under the $Spin^+(p,q)$-action as discussed in \cite{br}:

\begin{enumerate}
\item In signatures $(2,2)$ and $(3,3)$ every real half-spinor $\ph \in \Delta_{m,m}^{\R,\pm} \backslash \{0 \}$ is pure. Consequently, every real twistor half-spinor on $(M^{2,2},g)$ without zeroes is \textit{double pure}, by which we mean that both $\ph$ and the associated spin tractor $\psi$ are pointwise pure.
\item In signature $(3,2)$ every nonzero real spinor is pure. In signature $(4,3)$ a real spinor $\psi \in \Delta^{\R}_{4,3}$ is pure iff it is nonzero and $\langle \psi, \psi \rangle_{\Delta_{4,3}^{\R}} = 0$. With the scalar product formula (\ref{zz}) one concludes that a twistor spinor $\ph \in \Gamma(M^{3,2},S^g)$ is double pure iff $\ph$ admits no zeroes and $\langle \ph , D^g \ph \rangle \equiv 0$.
\item Similarly one shows that a twistor spinor $\ph \in \Gamma(M^{3,3},S^g)$ without zeroes is double pure iff $\langle \ph , D^g \ph \rangle \equiv 0$.
\end{enumerate}
Using Proposition \ref{tms} then directly yields the following:

\begin{Proposition} \label{propos}
Real twistor half-spinors in signature $(2,2)$ without zeroes and real twistor (half-)spinors without zeroes in signatures $(3,2)$ and $(3,3)$ satisfying that $\langle \ph , D^g \ph \rangle \equiv 0$ are locally conformally equivalent to parallel spinors (off a singular set). Their associated distributions ker $\ph \subset TM$ are integrable (off a singular set).
\end{Proposition}
Moreover, in the mentioned cases, the locally parallel spinor is real and pure at every point and the considered signatures are split signatures. In view of this, (\ref{purre1}) gives a local normal form for the metric.
Consequently, one has a complete local description of the geometries admitting \textit{non-generic} twistor spinors in signatures $(3,2)$ and $(3,3)$. This complements the classification of geometries admitting generic twistor spinors in signatures $(3,2)$ and $(3,3)$ from \cite{hs1}, i.e. $\langle \ph , D^g \ph \rangle \neq 0$, where the associated distribution ker $\ph$ turns out to be generic. Moreover, since the mentioned non-generic twistor spinors are pure, Proposition \ref{gogo} applies, yielding the following about the zero set structure:

\begin{Proposition}
Let $\ph \in \Gamma(S^g)$ be a real twistor (half-)spinor in signature $(2,2)$ or $(3,2)$. Then the zero set $Z_{\ph}$ -if nonempty- consists locally of totally lightlike planes. For a real twistor half-spinor with zero in signature $(3,3)$ the zero set is locally an embedded 3-dimensional totally lightlike submanifold.
\end{Proposition}

\subsection{Twistor spinors in signature (4,3)}
We start with algebraic observations: As known from \cite{kath}, $Spin^+(4,3)$ acts transitive on each of the level sets $M_c:=\{v \in \Delta_{4,3}^{\R} \mid \langle v, v \rangle = c \}$ for $c\neq 0$ and $M_0 \backslash \{0\}$ is precisely the space of pure spinors. Moreover, it holds for $\tau \in M_{c \neq 0}$ that ker $\tau=\{0 \}$ and $\langle X \cdot \tau, \tau \rangle = 0$ for all $X \in \R^{4,3}$. We use this in order to describe the orbit structure of $\Delta_{5,4}^{\R}$. According to \cite{br} each $N_c:=\{v \in \Delta_{5,4}^{\R} \mid \langle v, v \rangle = c \}$ constitutes for $c\neq 0$ a single orbit. However, $N_c \backslash \{0\}$ decomposes into at least 2 orbits wrt. the $Spin^+(5,4)-$action. 

\begin{Lemma}
It holds that ker $v \neq \{0\}$ for all $v \in N_{0} \subset \Delta_{5,4}^{\R}$.
\end{Lemma}

\begin{pf}We realise $\Delta_{5,4} = Ann(e_-) \oplus Ann(e_+) \cong \Delta_{4,3} \oplus \Delta_{4,3}$ as described in (\ref{fs}). With respect to this identification, we write $v = \begin{pmatrix} \tau \\ \chi \end{pmatrix}$ for $\tau, \chi \in \Delta_{4,3}^{\R}$. It follows that a vector $x=\alpha e_- + y+\beta e_+$ acts as
\begin{align}
x \cdot v = \begin{pmatrix} y\cdot \tau + \alpha \chi \\ -y \cdot \chi - 2 \beta \tau \end{pmatrix}. \label{g}
\end{align}
The scalar product formula (\ref{zz}) implies that $v \in N_{0} \Leftrightarrow \langle \tau, \chi \rangle_{\Delta_{4,3}^{\R}}=0$. If one of $\tau, \chi$ is trivial, the claim is obvious. Otherwise, we distinguish two cases:
Suppose that $\langle \chi, \chi \rangle = \langle \tau, \tau \rangle = 0$. In this case $\tau$ and $\chi$ are pure spinors with trivial pairing. It is a classical fact (cf. \cite{pure}) that in this case ker $\tau \cap \text{ ker }\chi \neq \{0\}$. (\ref{g}) implies that each nonzero element of this intersection lies in ker $v$. Thus it remains to consider the case where (wlog.) $\langle \tau, \tau \rangle \neq 0$. If $\chi \notin \R^{4,3}\cdot \tau$, we would for dimensional reasons have that $\Delta^{\R}_{4,3} = \R^{4,3}\cdot \tau \oplus \R\chi$ which implies that $\langle \Delta^{\R}_{4,3}, \tau \rangle = 0$ in contradiction to $\tau \neq 0$. Therefore, we find $y \in \R^{4,3}$ with $y \cdot \tau = \chi$. It follows that $y \cdot \chi + ||y||^2 \tau = 0$ yielding that $e_- -y -\frac{||y||^2}{2}e_+ \in \text{ker }v$. 
\end{pf}
Moreover, one calculates using orbit representatives that ker $v=0$ if $v \in N_{c \neq 0}$.
\begin{Proposition} \label{kj}
Let $\ph$ be a real twistor spinor on a conformal space $(M^{4,3},c)$. Let $g \in c$. Then exactly one of the following cases occurs:
\begin{enumerate}
\item It is $\langle \ph, D^g \ph \rangle \equiv 0$. In this case the spinor is locally equivalent to a parallel spinor off a singular set. One either has locally a parallel pure spinor field with a normal form of the metric given by \ref{purre1} or the spinor is locally a parallel spinor on a space whose holonomy representation is contained in $G_2$.
\item It is $\langle \ph, D^g \ph \rangle \neq 0$. Up to singular points there is locally around each point an Einstein metric with nonzero scalar curvature in the conformal class. The twistor spinor cannot be rescaled to a parallel spinor but decomposes into the sum of two Killing spinors.
\end{enumerate}
\end{Proposition}

\begin{pf}Let $\psi \in \Gamma(\mathcal{S})$ be the parallel spin tractor associated to $\ph$. The scalar product formula (\ref{zz}) shows that $\langle \psi, \psi \rangle =  \text{const.} \cdot \langle \ph, D^g \ph \rangle$. Thus, the previous algebraic discussion yields that $\langle \ph, D^g \ph \rangle \equiv 0$ iff $\psi$ admits a nontrivial kernel $\mathcal{H}_{\psi}$. It follows that $\ph$ is locally parallel by section 3. Thus, $\ph$ is locally of constant $Spin^+(4,3)$ orbit type. This means it is either pure or has trivial kernel. The stabilizer of a spinor $v \in M_{c \neq 0} \subset \Delta_{4,3}^{\R}$ is isomorphic to a copy of the exceptional group $G_2$ (cf. \cite{kath}). This proves the first part.\\
$\langle \ph, D^g \ph \rangle \neq 0$ is equivalent to $\langle \psi, \psi \rangle \neq 0$. It is known from \cite{br} that for every $v \in N_{c \neq 0 } \subset \Delta_{5,4}^{\R}$ one has
 \begin{align} \lambda(Stab_v Spin^+(5,4)) \cong Spin^+(4,3) \subset SO+(4,4) \subset SO+(5,4). \label{st}
\end{align}The conformal holonomy representation thus stabilizes a non-null vector yielding an Einstein metric on an open, dense subset. It is a classical fact (cf. \cite{boh}) that on an Einstein space every twistor spinor decomposes into the sum of two Killing spinors. 
$\hfill \Box$\\
\newline
The occuring local geometries admitting parallel spinors in signature $(4,3)$ are well understood (cf. \cite{kath}). Moreover, since dim ker $v \in \{0,3\}$ for all $v \in \Delta_{4,3} \backslash\{0\}$ one has that the zero set of a real twistor spinor in signature $(4,3)$ with zero consist locally either of isolated points or of 3-dimensional totally lightlike planes. In the first case one has locally $G_2$-holonomy off the zero set, in the second case one locally has a parallel pure spinor off the zero set as follows from the dimension of $L_{\ph}$ and the proof of the last Proposition.
\end{pf}

\subsection{Twistor spinors in signature (4,2)}
We again start with some algebraic observations. The complex spinor module $\Delta_{4,2}^{\C}$ admits a real structure commuting with Clifford multiplication giving the real module $\Delta_{4,2}^{\R}$\footnote{Note that in contrast to the complex case, $\Delta_{4,2}^{\R}$ is irreducible as $Spin^+(4,2)-$module. Hence there are no real half-spinors in signature $(4,2)$.}. The same applies to $\Delta_{5,3}^{\R}$. We consider the map
\begin{align*}
i: {\Delta_{5,4}^{\R}}_{|Spin^+(5,3)} \rightarrow \Delta_{5,3}^{\R}, \text{  }v \mapsto v
\end{align*}
being an isomorphism of $Spin^+(5,3)$-representations. In this way we can view $\Delta_{5,4}^{\R}$ as $Spin^+(5,3)-$ module and it holds that $\langle i(v),i(v) \rangle_{\Delta_{5,3}^{\R}}= \langle v,v \rangle_{\Delta_{5,4}^{\R}}$. Let $v \in \Delta_{5,3}^{\R}$ with $\langle v,v \rangle \neq 0$. (\ref{st}) then yields that
\begin{align*}
\lambda(Stab_v Spin^+(5,3)) \subset \lambda(Stab_v Spin^+(5,4)) \cong Spin^+(4,3) \subset SO^+(4,4) \subset SO^+(5,4). 
\end{align*}
However, as also $\lambda(Stab_v Spin^+(5,3)) \subset SO^+(5,3)$ we see that in fact up to conjugation \\
$\lambda(Stab_v Spin^+(5,3)) \subset SO^+(4,3) \subset SO^+(5,3)$. Thus there is a stabilized non-null vector in $\R^{5,3}$. \\
If $\langle v,v \rangle = 0$ we cannot make a general statement about ker $v$. There is the subcase of pure spinors but it is also possible for $v$ to have trivial kernel. There is no complete orbit classification available. \\
\newline
In complete analogy to the second case of Proposition \ref{kj} one now shows the following:
\begin{Proposition}
Let $\ph \in \Gamma(M^{4,2},S_{\R}^g)$ be a twistor spinor with $\langle \ph, D^g \ph \rangle \neq 0$. Then there is on an open, dense subset an Einstein metric with nonzero scalar curvature in the conformal class. Moreover, the spinor cannot be resclaed to a parallel spinor.
\end{Proposition}

We cannot completely describe the case $\langle \ph, D^g \ph \rangle = 0$. There is a subcase in which $\langle \ph, \ph \rangle \equiv 0 \equiv \langle D^g \ph , D^g \ph \rangle $ in which $\ph$ is locally conformally equivalent to a parallel spinor. This follows since the assumptions gurantee the existence of $0 \neq X \in TM$ with $X \cdot \ph = X \cdot D^g \ph = 0$ which implies that $\mathcal{H}_{\psi} \neq 0$. There is another subcase when  $\ph$ is parallel and $L_{\ph}=\{0\}$ where one has a Ricci-flat pseudo-Kaehler metric in the conformal class (cf. the discussion of the $p=2$ case at the end of section 3). \\
For a nonzero spinor $v \in \Delta_{4,2}^{\R}$ it holds that dim ker $v \in \{0,2\}$. Thus, the zero set of a real twistor spinor with zeroes in signature $(4,2)$ consists either of totally lightlike nullplanes, where one has locally a parallel pure spinor off the zero set, or $Z_{\ph}$ consists of isolated points and the geometry off the zero set is Ricci-flat pseudo-Kaehler.

\section*{Acknowledgement}
The author gladly acknowledges support from the DFG (SFB 647 - Space Time Matter at Humboldt University Berlin) and the DAAD (German Academic Exchange Service).

\bibliographystyle{elsarticle-num}
\bibliography{literatur}
\end{document}